\documentclass[MSNbibl,number,citesort,seceqn,dvips]{arxbj}
\usepackage{upgreek}

\aid{0}
\volume{19}
\issue{5B}
\pubyear{2013}
\firstpage{2167}
\lastpage{2199}
\doi{10.3150/12-BEJ448} 

\makeatletter

\newcommand{\rrvert}{\vert}
\newcommand{\llvert}{\vert}

\newtheorem{tm}{Theorem}[section]
\newremark{remark}[tm]{Remark}
\newproclaim{definition}[tm]{Definition}
\newtheorem{corollary}[tm]{Corollary}
\newtheorem{proposition}[tm]{Proposition}
\newtheorem{lemma}[tm]{Lemma}

\newcommand{\sgn}{\operatorname{sgn}}
\renewcommand{\Re}{\operatorname{Re}}

\newcommand{\R}{\mathbb{R}}
\newcommand{\ZZ}{\mathbb{Z}}
\newcommand{\N}{\mathbb{N}}
\newcommand{\CC}{\mathbb{C}}
\makeatother

\begin{document}
\begin{frontmatter}

\title{Recurrence and transience property for a class of Markov chains}
\runtitle{Recurrence and transience property for a class of Markov chains}

\begin{aug}
\author{\fnms{Nikola} \snm{Sandri\'{c}}\corref{}\ead[label=e1]{nsandric@grad.hr}}
\runauthor{N. Sandri\'{c}} 
\address{Department of Mathematics,
Faculty of Civil Engineering, University of Zagreb, Fra Andrije Ka\v
ci\'ca-Mio\v si\'ca 26, 10 000 Zagreb, Croatia. \printead{e1}}
\end{aug}

\received{\smonth{3} \syear{2011}}
\revised{\smonth{4} \syear{2012}}

%
\begin{abstract}
We consider the recurrence and transience problem for a
time-homogeneous Markov chain on the real line with transition
kernel $p(x,\mathrm{d}y)=f_x(y-x)\, \mathrm{d}y$, where the density
functions $f_x(y)$,
for large $|y|$, have a power-law decay with exponent $\alpha(x)+1$,
where $\alpha(x)\in(0,2)$.
In this paper, under a uniformity condition on the density functions
$f_x(y)$ and an additional mild drift
condition, we prove that when
$\lim\inf_{|x|\longrightarrow\infty}\alpha(x)>1$, the chain is
recurrent. Similarly, under the same uniformity condition on the
density functions $f_x(y)$ and some mild technical conditions, we
prove that when $\lim\sup_{|x|\longrightarrow\infty}\alpha(x)<1$, the
chain is transient. As a special case of these results, we give a new
proof for the recurrence and transience property of a symmetric
$\alpha$-stable random walk on $\R$ with the index of stability
$\alpha\in(0,1)\cup(1,2).$
\end{abstract}

%

\end{frontmatter}

\section{Introduction}

Let $(\Omega,\mathcal{F},\mathbb{P})$ be a probability space and let
$\{Z_n\}_{n\in\N}$ be a sequence of i.i.d. random variables on
$(\Omega,\mathcal{F},\mathbb{P})$ taking values in $\R^{d}$,
$d\geq1$. Let us define $X_n:=\sum_{i=1}^{n}Z_i$ and $X_0:=0$. The
sequence $\{X_n\}_{n\geq0}$ is called a \emph{random walk} with
jumps $\{Z_n\}_{n\in\N}$. The random walk $\{X_n\}_{n\geq0}$ is
said to be \emph{recurrent} if
\[
\mathbb{P} \Bigl(\operatorname{\lim\inf}\limits_{n\longrightarrow\infty}|X_n|=0 \Bigr)=1,
\]
and \emph{transient} if
\[
\mathbb{P} \Bigl(\lim_{n\longrightarrow\infty}|X_n|=\infty \Bigr)=1.
\]
It is well known that every random walk is either recurrent or
transient (see \cite{durrett}, Theorem 4.2.1). Recall that a random
walk $\{X_n\}_{n\geq0}$ in $\R^{d}$ is called \emph{truly
$d$-dimensional} if $\mathbb{P}(\langle Z_1, x\rangle\neq0)>0$ holds
for all $x\in\R^{d}\setminus\{0\}.$ It is also well known that every
truly $d$-dimensional random walk is transient if $d\geq3$
(see \cite{durrett}, Theorem 4.2.13). An $\R^{d}$-valued random
variable $Z$ is said to have \emph{stable distribution} if, for any
$n\in\N$, there are $a_n
> 0$ and $b_n\in\R^{d}$, such that
\[
Z_1+\cdots+Z_n\stackrel{\mathrm{d}} {=}a_nZ+b_n,
\]
where
$Z_1,\ldots, Z_n$ are independent copies of $Z$ and
$\stackrel{\mathrm{d}}{=}$ denotes equality in
distribution. It turns out that $a_n=n^{{1}/{\alpha}}$ for some
$\alpha\in(0,2]$ which is called the index of stability
(see \cite{taqqu}, Definition~1.1.4 and Corollary 2.1.3). The case
$\alpha=2$ corresponds to the Gaussian random variable. A random
walk $\{X_n\}_{n\geq0}$ is said to be stable if the random variable
$Z_1$ has stable distribution. In the class of truly two-dimensional
stable random walks in $\R^{2}$, by \cite{durrett}, Theorem
4.2.9, the only recurrent case is the case when
$\{X_n\}_{n\geq0}$ is a truly two-dimensional random walk with zero
mean Gaussian jumps. In the case $d=1$, every stable distribution is
characterized by four parameters: the stability parameter
$\alpha\in(0,2],$ the skewness parameter $\beta\in[-1,1]$, the
scale parameter $\gamma\in(0,\infty)$ and the shift parameter
$\delta\in\R$ (see \cite{taqqu}, Definition 1.1.6). Using the
notation from \cite{taqqu}, we denote one-dimensional stable
distributions by $S_{\alpha}(\beta,\gamma,\delta).$ For symmetric
stable distributions, that is, for $S_{\alpha}(0,\gamma,0)$
(see \cite{taqqu}, Property 1.2.5), we write S$\alpha$S. A
S$\alpha$S random walk is recurrent if and only if $\alpha\geq1$
(see the discussion after \cite{durrett}, Lemma 4.2.12). In this
paper, we generalize the S$\alpha$S random walk in the way that the
index of stability of the jump distribution depends on the current
position and study the transience and recurrence property of the
generalization.

Actually, we will not need the stability property of transition
jumps. All we will need is a tail behavior of transition jumps. Let
us introduce the notation $f(y)\sim g(y),$ when $y\longrightarrow
y_0$, for $\lim_{y\longrightarrow y_0}f(y)/g(y)=1,$ where
$y_0\in[-\infty,\infty]$. Recall that if $f(y)$ is the density
function of a S$\alpha$S distribution with $\alpha\in(0,2)$ and
$\gamma\in(0,\infty)$ (for the existence of densities of
$S_{\alpha}(\beta,\gamma,\delta)$ distributions see
\cite{taqqu}, Definition 1.1.6 and \cite{durrett}, Theorem 3.3.5),
then
\[
f(y)\sim c_{\alpha}|y|^{-\alpha-1},
\]
when
$|y|\longrightarrow\infty,$ where $c_1=\frac{\gamma}{2}$ and
$c_{\alpha}=\frac{\gamma}{\uppi}\Gamma(\alpha+1)\sin (\frac
{\uppi\alpha
}{2} ),$
for $\alpha\neq1,$ see \cite{taqqu}, Property~1.2.15.
Now, let
$\alpha\dvtx \R\longrightarrow(0,2)$ and $c\dvtx \R\longrightarrow(0,\infty)$
be arbitrary functions and let $\{f_x\dvt x\in\R\}$ be a family of
density functions on $\R$ such that
\begin{enumerate}[(C2)]
\item[(C1)] $x\longmapsto f_x(y)$ is a Borel measurable
function for
all $y\in\R$ and
\item[(C2)] $f_x(y)\sim
c(x)|y|^{-\alpha(x)-1},$ when $|y|\longrightarrow\infty,$ for all
$x\in\R$.
\end{enumerate}
Let us define a Markov chain $\{X_n\}_{n\geq0}$ on $\R$ by the
following transition kernel
%
\begin{equation}\label{eq11}
p(x,\mathrm{d}y):=f_x(y-x)\,\mathrm{d}y.
\end{equation}
The chain $\{X_n\}_{n\geq0}$ jumps from the
state $x$ with transition density $f_x(y-x)$, with the power-law
decay with exponent $\alpha(x)+1$, and this jump distribution
depends only on the current state $x$. Transition densities
$\{f_x\dvt x\in\R\}$ are asymptotically equivalent to the densities of
S$\alpha$S distributions, and we call such chain
a \emph{stable-like chain}. The aim of this paper is to find
conditions for the recurrence and
transience property of the stable-like chain $\{X_n\}_{n\geq0}$ in
terms of the function $\alpha(x).$

To the best of our knowledge, all methods used in establishing
conditions for recurrence and transience in the random walk case
are based on the i.i.d. property of random walk jumps, that is, laws of
large numbers (Chung--Fuchs theorem), central limit theorems,
characteristic functions approach (Stone--Ornstein formula) etc.
(see \cite{durrett}, Theorems 4.2.7, 4.2.8 and 4.2.9). Although we
deal with distributions similar to S$\alpha$S distributions, it is
not clear if these methods can be used in the case of the
non-constant function $\alpha(x)$.

Special cases of
this problem have been considered in \cite{bjoern,franke,frankeerata,kemperman} and \cite{rogozin}.
In \cite{kemperman} and \cite{rogozin}, the authors consider
the countable state space $\ZZ$ and the function $\alpha(x)$ is a
two-valued step function which takes one value on negative integers
and the other one on nonnegative integers. The processes considered
in \cite{bjoern} and \cite{franke} run in continuous time. The
function $\alpha(x)$ considered in \cite{bjoern} is a two-valued
step function which takes one value on negative reals and the other
one on nonnegative reals, while in \cite{franke} the author
considers the case when the function $\alpha(x)$ is periodic and
continuously differentiable. The methods used in \cite{kemperman,rogozin,bjoern} and \cite{franke}, actually reduce the
process to random walks and L\'{e}vy processes. Also, it is not
clear if these methods can be used in the general case, that is, when
the function $\alpha(x)$ is an arbitrary function. In this paper,
under certain assumptions on the functions $\alpha(x)$, $c(x)$ and
on the family of density functions $\{f_x\dvt x\in\R\}$, we give
sufficient conditions for the recurrence and transience property of
the stable-like chain $\{X_n\}_{n\geq0}$ in terms of the function
$\alpha(x).$

Let us denote by $\mathcal{B}(\R)$ the Borel $\sigma$-algebra on
$\R$, by $\lambda$ the Lebesgue measure on $\mathcal{B}(\R)$ and for
arbitrary $B\in\mathcal{B}(\R)$ and $x\in\R$ we define
$B-x:=\{y-x\dvt y\in B\}$. Assume that the family of probability
densities $\{f_x\dvt x\in\R\}$ satisfies additional three conditions:
\begin{enumerate}[(C4)]
\item[(C3)] there exists $k>0$ such that
\[
\lim_{|y|\longrightarrow\infty}\sup_{x\in[-k,k]^{c}}\biggl\llvert f_x(y)
\frac
{|y|^{\alpha(x)+1}}{c(x)}-1\biggr\rrvert =0;
\]
\item[(C4)] $\inf_{x\in C}c(x)>0$
for every compact set
$C\subseteq[-k,k]^{c}$;
\item[(C5)] there exists $l>0$ such that for every compact
set $C\subseteq[-l,l]^{c}$ with $\lambda(C)>0$,
we have
\[
\inf_{x\in[-k,k]}\int_{C-x}f_x(y)\,\mathrm{d}y>0.
\]
\end{enumerate}

Condition (C3) ensures that out of some compact set all jump
densities of the stable-like chain $\{X_n\}_{n\geq0}$ can be
replaced by their tail behavior uniformly. This condition is crucial
in proving certain structural properties of the chain
$\{X_n\}_{n\geq0}$ and in finding sufficient conditions for the
recurrence and transience. Another essential property of the chain
$\{X_n\}_{n\geq0}$ is that every compact set is a petite set. A
petite set is a set which assumes a role of a singleton for Markov
chains on general state space (for the exact definition of the
petite set see Definition \ref{d22}). This is the reason why
compact sets are important in conditions (C3), (C4) and (C5).
Besides ensuring that all compact sets are
petite sets (singletons), conditions (C4) and (C5) ensure also that
the chain is irreducible. Condition (C4) ensures that the scaling
function $c(x)$ does not vanish on petite sets, and condition (C5)
ensures that the petite set $[-k,k]$ communicates with the rest of
the state space.

\begin{remark} Note that condition (C3) implies
%
\begin{equation}
\label{eq12}\sup_{x\in[-k,k]^{c}}c(x)<\infty.
\end{equation}
Indeed, let
$0<\varepsilon<1$ be arbitrary. Then there exists
$y_{\varepsilon}\geq1$ such that for all $|y|\geq y_{\varepsilon}$
we have
\[
\biggl\llvert f_x(y)\frac{|y|^{\alpha(x)+1}}{c(x)}-1\biggr\rrvert <\varepsilon
\]
for all $x\in[-k,k]^{c}$. Therefore, upon integrating over $y$ we
get
\[
c(x)<\frac{1}{1-\varepsilon} \biggl(2\int_{y_{\varepsilon}}^{\infty
}y^{-\alpha(x)-1}\,\mathrm{d}y
\biggr)^{-1}\leq\frac{1}{1-\varepsilon} \biggl(2\int_{y_{\varepsilon}}^{\infty}y^{-3}\,\mathrm{d}y
\biggr)^{-1}=\frac
{y^{2}_{\varepsilon
}}{1-\varepsilon}
\]
for every
$x\in[-k,k]^{c}$.
\end{remark}
An example of a stable-like chain which satisfies conditions
(C3)--(C5) is the chain which has exactly
$S_{\alpha(x)}(0,\gamma(x),\delta(x))$ jumps at each location $x$,
where the functions $\alpha(x)$, $\gamma(x)$ and $\delta(x)$ are
Borel measurable and take finitely many values (see Proposition
\ref{p55} for details).

Before stating the main results of this paper we recall relevant
definitions of recurrence and transience.
%
\begin{definition}\label{d12} Let $\{Y_n\}_{n\geq0}$ be a Markov
chain on
$(\R,\mathcal{B}(\R))$.
\begin{enumerate}[(iii)]
\item[(i)] The chain $\{Y_n\}_{n\geq0}$ is \emph{$\varphi
$-irreducible} if there
exists a probability measure $\varphi$ on $\mathcal{B}(\R)$ such
that for every $x\in\R$ there exists $n\in\N$ such that
$\varphi(B)>0$ implies $\mathbb{P}(Y_n\in B|Y_0=x)>0$.
\item[(ii)] The chain $\{Y_n\}_{n\geq0}$ is \emph{recurrent} if it is
$\varphi$-irreducible and if $\sum_{n=0}^{\infty}\mathbb{P}(Y_n\in
B|Y_0=x)=\infty$ holds for all $x\in
\R$ and all $B\in\mathcal{B}(\R)$, such that $\varphi(B)>0$.
\item[(iii)] The chain $\{Y_n\}_{n\geq0}$ is \emph{transient} if it is
$\varphi$-irreducible and
if there exists a countable cover of $\R$ with sets
$\{B_j\}_{j\in\N}\subseteq\mathcal{B}(\R)$, such that for each
$j\in\N$ there is a finite constant $M_j\geq0$ such that
$\sum_{n=0}^{\infty}\mathbb{P}(Y_n\in B_j|Y_0=x)\leq M_j$ holds for
all $x\in\R$.
\end{enumerate}
\end{definition}
The following two constants will appear in the statements of the
main results: For $\alpha\in(1,2)$, let
\[
R(\alpha):=\sum_{i=1}^{\infty}
\frac{1}{i(2i-\alpha)}-\frac
{\ln
2}{\alpha}-\frac{1}{2\alpha} \biggl(\Psi \biggl(
\frac{\alpha
+1}{2} \biggr)-\Psi \biggl(\frac{\alpha}{2} \biggr) \biggr),
\]
and for $\alpha\in[0,1)$ and $\beta\in(0,1-\alpha)$ let
\[
T(\alpha,\beta):={}_2F_1(-\alpha,\beta,1-\alpha;1)+\beta
B(1;\alpha+\beta,1-\alpha)-\alpha B(1;\alpha+\beta,1-\beta),
\]
where
$\Psi(z)$ is the Digamma function, ${}_2F_1(a,b,c;z)$ is the Gauss
hypergeometric function and $B(x;z,w)$ is the incomplete Beta
function (see Section \ref{sec3} for the definition of these functions). The
constants $R(\alpha)$ and $T(\alpha,\beta)$ are strictly positive
(see proofs of Theorems \ref{tm13} and \ref{tm14}). Furthermore,
it is not hard to see that the constant $R(\alpha)$, as a function
of $\alpha\in(1,2)$, is strictly increasing, $R(1)=0$ and
$\lim_{\alpha\longrightarrow2}R(\alpha)=\infty.$ The constant
$T(\alpha,\beta)$, as a function of $\beta\in(0,1-\alpha)$ for fixed
$\alpha\in(0,1)$, is strictly positive and
$T(\alpha,0)=T(\alpha,1-\alpha)=0,$ while considered as a function
of $\alpha\in[0,1-\beta)$ for fixed $\beta\in(0,1)$, it is strictly
decreasing, $T(0,\beta)=2$ and $T(1-\beta,\beta)=0.$

\begin{tm}\label{tm13}
Let $\alpha\dvtx \R\longrightarrow(1,2)$ be an arbitrary function such
that
\[
\alpha:=\operatorname{\lim\inf}\limits_{|x|\longrightarrow\infty}\alpha(x)>1.
\]
Furthermore, let $c\dvtx \R\longrightarrow(0,\infty)$ be an arbitrary
function and let $\{f_x\dvt x\in\R\}$ be a family of density functions
on $\R$ which satisfies conditions \emph{(C1)--(C5)} and
such that
%
\begin{equation}
\label{eq13}\operatorname{\lim\sup}\limits_{|x|\longrightarrow\infty}\sgn(x)\frac{|{x}|^{\alpha({x})-\mathrm{1}}}{c(x)}
\mathbb{E}[X_{\mathrm{1}}-X_{\mathrm{0}}|X_{\mathrm{0}}=x]<R(\alpha)
\end{equation}
when $\alpha<2$, and the left-hand side in (\ref{eq13}) is
finite when $\alpha=2$. Then the stable-like Markov chain
$\{X_n\}_{n\geq0}$ given by the transition kernel
\[
p(x,\mathrm{d}y)=f_x(y-x)\,\mathrm{d}y,
\]
is recurrent.
\end{tm}

\begin{tm}\label{tm14}
Let $\alpha\dvtx \R\longrightarrow(0,1)$ be an arbitrary function such
that
\[
\alpha:=\operatorname{\lim\sup}\limits_{|x|\longrightarrow\infty}\alpha(x)<1
\]
and let $\beta\in(0,1-\alpha)$ be arbitrary. Furthermore, let
$c\dvtx \R\longrightarrow(0,\infty)$ be an arbitrary function and let
$\{f_x\dvt x\in\R\}$ be a family of density functions which satisfies
conditions \emph{(C1)--(C5)} and there exists $a_0>0$, such that
%
\begin{equation}
\label{eq14}\operatorname{\lim\inf}\limits_{|x|\longrightarrow\infty}\frac
{\alpha(x)|x|^{\alpha(x)}}{c(x)}\int_{-a}^{a}
\biggl(1- \biggl(1+\sgn (x)\frac{y}{{1}+|{x}|} \biggr)^{-\beta}
\biggr)f_{x}(y)\,\mathrm{d}y>-T(\alpha ,\beta)
\end{equation}
for all $a\geq a_0.$ Then the stable-like Markov chain
$\{X_n\}_{n\geq0}$ given by the transition kernel
\[
p(x,\mathrm{d}y)=f_x(y-x)\,\mathrm{d}y,
\]
is transient.
\end{tm}

Actually, instead of condition
(\ref{eq13}), in the proof of Theorem \ref{tm13}, we use the
following more technical but equivalent condition
%
\begin{equation}
\label{eq15}\operatorname{\lim\sup}\limits_{\delta\longrightarrow0}\operatorname{\lim\sup}
\limits_{|x|\longrightarrow\infty}\frac{(1+|x|)^{\alpha(x)}}{c(x)}\int
_{-\delta
(1+|x|)}^{\delta(1+|x|)}\ln \biggl(1+\sgn(x)
\frac{y}{1+|x|} \biggr)f_{x}(y)\,\mathrm{d}y<R(\alpha)
\end{equation}
(see
Section \ref{sec5} for details). Conditions (\ref{eq13}) (i.e.,
(\ref{eq15})) and (\ref{eq14}) are needed to control the
behavior of the family of density functions $\{f_x\dvt x\in\R\}$ on sets
symmetric around the origin. Condition (\ref{eq13}) actually says
that when the chain $\{X_n\}_{n\geq0}$ has moved far away from the
origin, since $R(\alpha)>0$, it cannot have strong tendency to move
further from the origin. Since $R(\alpha)>0$, it is clear that
condition (\ref{eq13}) is satisfied if $\alpha(x)\in(1,2)$ and if
$f_x(y)=f_x(-y)$ holds for all $y\in\R$ and for all $|x|$ large
enough. For a non-symmetric example, one can take $f_x(y)$ to be the
density function of a $S_{\alpha_-}(0,\gamma_-,\delta_-)$
distribution, when $x<0$, and the density function of a
$S_{\alpha_+}(0,\gamma_+,\delta_+)$ distribution, when $x\geq0$,
where $\alpha_-,\alpha_+\in(1,2)$, $\gamma_-,\gamma_+\in(0,\infty),$
$\delta_-\geq0$ and $\delta_+\leq0$.

Using the concavity property of the function $x\longmapsto
x^{\beta}$, for $\beta\in(0,1-\alpha)$, condition (\ref{eq14})
follows from the condition
%
\begin{equation}
\label{eq17}\operatorname{\lim\sup}\limits_{|x|\longrightarrow\infty}\frac
{\alpha
(x)}{c(x)}|x|^{\alpha(x)-1}<
\frac{T(\alpha,\beta)}{a_0\beta}
\end{equation}
(see Section \ref{sec5} for details). Note that condition (\ref{eq17})
actually says that the function $c(x)$ cannot decrease too fast.
Since $T(\alpha,\beta)>0$ and $\alpha(x)\in(0,1)$, a simple example
which satisfies condition (\ref{eq17}) is the case when
$c(x)\geq
d|x|^{\alpha(x)-1+\epsilon}$, for some $d>0$ and for all $|x|$ large
enough, where $0<\epsilon<1-\alpha$ is arbitrary. Furthermore, one
can prove that the function $\beta\longmapsto T(\alpha,\beta)/\beta$
is strictly decreasing on $(0,1-\alpha)$. Hence, according to the
condition (\ref{eq17}), we choose $\beta$ close to $0$.

In the random walk case, that is, when the family of density functions
$\{f_x\dvt x\in\R\}$ is reduced to a single density function $f(y)$ such
that $f(y)\sim c|y|^{-\alpha-1}$, when $|y|\longrightarrow\infty$,
where $\alpha\in(0,2)$ and $c\in(0,\infty)$,
conditions (C1)--(C5) are trivially satisfied. Hence, by Theorem \ref
{tm13} and the condition (\ref{eq13}), if $\alpha>1$ and if
\[
\int_{\R}y f(y)\,\mathrm{d}y=0,
\]
the random walk with the
jump density $f(y)$ is recurrent, and if $\alpha<1$, by
Theorem \ref{tm14} and the condition (\ref{eq17}), the random
walk with the jump density $f(y)$ is transient. This result can be
strengthened. If we assume that $f(y)=f(-y)$ for all $y\in\R$, from
the discussion in \cite{spitzer}, page 88, the random walk with the
jump density $f(y)$ is recurrent if and only if $\alpha\geq1.$ As a
simple consequence of Theorems \ref{tm13} and \ref{tm14}, we get
the following well-known recurrence and transience conditions for
the S$\alpha$S random walk case.
%
\begin{corollary}\label{c15}
A S$\alpha$S, $1<\alpha<2$, random walk is recurrent.
A $S_{\alpha}(0,\gamma,\delta)$, $0<\alpha<1$, random walk
with arbitrary shift is transient.
\end{corollary}
The previous corollary can be generalized. If the functions
$\alpha(x),$ $\gamma(x)$ and $\delta(x)$ are Borel measurable and
take finitely many values, then the stable-like chain with
S$\alpha(x)$S jumps is recurrent if $\alpha(x)\in(1,2)$
for all $x\in\R$. If
$\alpha(x)\in(0,1)$ for all $x\in\R$, then the stable-like chain
with $S_{\alpha(x)}(0,\gamma(x),\delta(x))$ jumps is transient.

\begin{remark}
Conditions in Theorems \ref{tm13} and \ref{tm14} are only sufficient
conditions for recurrence and transience
of the stable-like chain $\{X_n\}_{n\geq0}$. On the countable state
space $\ZZ$, when
\[
\alpha(i)=\cases{
\alpha, &\quad $i<0$,
\cr
\beta,&\quad $i\geq0$}
\]
for $\alpha,\beta\in(0,2)$, in \cite{kemperman,rogozin} it is
proved that if $\frac{\alpha+\beta}{2}>1$, the associated chain is
recurrent, and if $\frac{\alpha+\beta}{2}<1$, the associated chain
is transient. A similar result, with
\[
\alpha(x)=\cases{
\alpha, &\quad $x<0$,
\cr
\beta,&\quad $x\geq0$ }
\]
for $\alpha,\beta\in(0,2)$, is proved in the continuous time case in
\cite{bjoern}, that is, a stable-like process with the symbol
$|\xi|^{\alpha(x)}$ is recurrent if and only if
$\frac{\alpha+\beta}{2}\geq1.$ In \cite{franke}, in the case when
the function $\alpha(x)$ is periodic and continuously differentiable
function, it is proved that all that matters is the minimum of the
function $\alpha(x)$. If
$\lambda(\{x\dvt \alpha(x)=\alpha_0:=\inf\{\alpha(y)\dvt y\in\R\}\})>0$,
then a stable-like process with the symbol $|\xi|^{\alpha(x)}$ is
recurrent if and only if $\alpha_0\geq1.$
\end{remark}

Now we explain our strategy of proving the main results. The proof
of Theorems \ref{tm13} and \ref{tm14} is based on the
\emph{Foster--Lyapunov drift criterion} for recurrence and transience
of Markov chains (see \cite{meyn}, Theorems 8.4.2 and 8.4.3). This
criterion is based on finding an appropriate test function $V(x)$
(positive and unbounded in the recurrence case and positive and
bounded in the transience case), and an appropriate set
$C\in\mathcal{B}(\R)$ (petite set) such that
$\int_{\R}p(x,\mathrm{d}y)V(y)-V(x)\leq0$, in the recurrence case, and
$\int_{\R}p(x,\mathrm{d}y)V(y)-V(x)\geq0$, in the transience case, for every
$x\in C^{c}.$ The idea is to find test functions $V(x)$ such that
the associated level sets $C_V(r):=\{y\dvt V(y)\leq r\}$ are compact
sets, that is, petite sets, and that $C_V(r)\uparrow\R$, when
$r\longrightarrow\infty$, in the case of recurrence and
$C_V(r)\uparrow\R$, when $r\longrightarrow1$, in the case of
transience. In the recurrence case for the test function, we take
$V(x)=\ln(1+|x|)$, and in the transience case we take
$V(x)=1-(1+|x|)^{-\beta},$ where $0<\beta<1-\alpha$ (recall that
$\alpha=\operatorname{\lim\sup}_{|x|\longrightarrow\infty}\alpha(x)<1$). Now, by
proving
that
\[
\operatorname{\lim\sup}\limits_{|x|\longrightarrow\infty}\frac{|x|^{\alpha
(x)}}{c(x)} \biggl(\int_{\R}p(x,\mathrm{d}y)V(y)-V(x)
\biggr)<0,
\]
in the recurrence case, and
\[
\operatorname{\lim\inf}\limits_{|x|\longrightarrow\infty}\frac{\alpha(x)|x|^{\alpha
(x)+\beta
}}{c(x)} \biggl(\int_{\R}p(x,\mathrm{d}y)V(y)-V(x)
\biggr)>0,
\]
in the transience case, since compact sets are petite sets, the
proofs of Theorems \ref{tm13} and \ref{tm14} are accomplished.

A similar approach, by using similar test functions $V(x)$, can be
found in \cite{lamperti} and \cite{rus}. In \cite{lamperti}, the
author considers a Markov chain on the nonnegative real line with
uniformly bounded transition jumps, while in \cite{rus} the
authors generalize this result to the case of uniformly bounded
$2+\delta_0$-moments of transition jumps, for some $\delta_0>0$. If
we allow that $\alpha(x)\in(0,\infty)$ and assume the following
additional assumption: $\sup_{x\in C}\alpha(x)<\infty$, for every
compact set
$C\subseteq[-k,k]^{c}$ (recall that the constant $k$ is defined in condition
(C3)), one can prove all nice structural
properties of the chain $\{X_n\}_{n\geq0}$, given by (\ref{eq11}),
proved in Section \ref{sec2}.
Hence, since the chain $\{X_n\}_{n\geq0}$ is recurrent if and only
if the chain $\{|X_n|\}_{n\geq0}$ is recurrent, \cite{rus} covers
the case when $\operatorname{\lim\inf}_{|x|\longrightarrow\infty}\alpha(x)>2.$

The paper is organized as follows. In Section \ref{sec2}, we give several
structural properties of the stable-like chain $\{X_n\}_{n\geq0}$
which will be crucial in finding sufficient conditions for the
recurrence and transience property. In Sections \ref{sec3} and \ref{sec4}, using
Foster--Lyapunov drift criterion for recurrence and transience of
Markov chains, we prove Theorems \ref{tm13} and \ref{tm14}. In
Section \ref{sec5}, we extend our model from the model of asymptotically
symmetric transition jumps to the model of asymptotically
non-symmetric transition jumps. Further, we prove that the change of
the chain $\{X_n\}_{n\geq0}$ on bounded sets will not affect the
recurrence and transience property.

Throughout the paper, we use the following notation. We write
$\ZZ_+$ and $\ZZ_-$ for nonnegative and nonpositive integers,
respectively. For $x,y\in\R$ let $x\wedge y =\min\{x,y\}$ and $x\vee
y =\max\{x,y\}$. Furthermore, $\{X_n\}_{n\geq0}$ will denote the
stable-like Markov chain on $\R$ given by (\ref{eq11}) with
transition densities satisfying conditions (C1)--(C5), while
$\{Y_n\}_{n\geq0}$ will denote an arbitrary Markov chain on
$(\R,\mathcal{B}(\R))$ given by the transition kernel $p(x,B),$ for
$x\in\R$ and $B\in\mathcal{B}(\R)$. For $x\in\R$, $B\in
\mathcal{B}(\R)$ and $n\in\N$ let $p^{n}(x,B):=\mathbb{P}(Y_n\in
B|Y_0=x)$ and $\tau_B:=\min\{n\geq1\dvt Y_n\in B\}$.

\section{Structural properties of the model}\label{sec2}

In this section, we discuss several structural properties of
stable-like Markov chains. In Definition \ref{d12}, we defined
irreducibility of a Markov chain on the state space
$(\R,\mathcal{B}(\R))$. In \cite{meyn}, Proposition~4.2.1, it is
shown that the irreducibility measure can always be maximized, that
is, if $\{Y_n\}_{n\geq0}$ is a $\varphi$-irreducible Markov chain,
then there exists
a probability measure $\psi$ on $\mathcal{B}(\R)$ such that the chain
$\{Y_n\}_{n\geq0}$ is
$\psi$-irreducible and $\varphi'\ll\psi$, for every irreducibility
measure $\varphi'$ on $\mathcal{B}(\R)$ of the chain
$\{Y_n\}_{n\geq0}$. The measure $\psi$ is called the \emph{maximal
irreducibility measure} and from now on, when we refer to
irreducibility measure we actually refer to the maximal
irreducibility measure. For the $\psi$-irreducible Markov chain
$\{Y_n\}_{n\geq0}$ on $(\R,\mathcal{B}(\R))$, let us set
$\mathcal{B}^{+}(\R)=\{B\in\mathcal{B}(\R)\dvt \psi(B)>0\}.$

\begin{proposition}\label{p21} Under conditions \emph{(C1)--(C4)}, the
maximal irreducibility measure for the chain $\{X_n\}_{n\geq0}$ is equivalent,
in the absolutely continuous sense, with the Lebesgue measure.
Therefore, the chain $\{X_n\}_{n\geq0}$ is $\lambda$-irreducible.
\end{proposition}
\begin{pf}
First, we prove that under conditions (C1)--(C4), the chain
$\{X_n\}_{n\geq0}$ is $\varphi$-irreducible for all measures
$\varphi$, such that $\varphi\ll\lambda$. We prove that for every
$x\in\R$ and for every $B\in\mathcal{B}(\R)$, such that
$\lambda(B)>0$, there exists $n\in\N$, such that $p^{n}(x,B)>0.$ It
is enough to prove the claim in the case of bounded sets. Let $B\in
\mathcal{B}(\R)$, $\lambda(B)>0,$ be an arbitrary bounded set. Let
$x\in\R$ and $0<\varepsilon<1$ be arbitrary. Then, by (C2), there
exists $y_{\varepsilon,x}\geq1$ such that for all $|y|\geq
y_{\varepsilon,x}$ we have
\[
\biggl\llvert f_x(y)\frac{|y|^{\alpha(x)+1}}{c(x)}-1\biggr\rrvert <\varepsilon.
\]
Furthermore, by (C3), there exists $k>0$ such that for given
$\varepsilon$ there exists $y_{\varepsilon}\geq1$, such that for all
$|y|\geq y_{\varepsilon}$ and all $z\in[-k,k]^{c}$, we have
\[
\biggl\llvert f_z(y)\frac{|y|^{\alpha(z)+1}}{c(z)}-1\biggr\rrvert <\varepsilon.
\]
Let $a:=\sup B$ and $y_{0}:=(y_{\varepsilon,x}\vee
y_{\varepsilon}\vee k) +|x|+|a|+1$. Finally, by (C4) we have
\begin{eqnarray*}
p^{2}(x,B)&=&\int_{\R}p(x,\mathrm{d}y)p(y,B)=\int
_{\R
}f_x(y-x)\int_{B-y}f_y(z)\,\mathrm{d}z\,\mathrm{d}y
\\
&\geq&\int_{y_{0}}^{2y_0}f_x(y-x)\int
_{B-y}f_y(z)\,\mathrm{d}z\,\mathrm{d}y
\\
&>&(1-\varepsilon)^{2}c(x)\int_{y_{0}}^{2y_{0}}(y-x)^{-\alpha
(x)-1}c(y)
\int_{B-y}|z|^{-\alpha(y)-1}\,\mathrm{d}z\,\mathrm{d}y
\\
&>&(1-\varepsilon)^{2}c(x) \Bigl(\inf_{y_{0}\leq y\leq2
y_{0}}c(y) \Bigr)\int
_{y_{0}}^{2y_{0}}(y-x)^{-3}\int
_{B-y}|z|^{-3}\,\mathrm{d}z\,\mathrm{d}y>0,
\end{eqnarray*}
since $B-y\subseteq(-\infty,-y_{\varepsilon})$, for $y\geq y_{0}$.

Now, we show the maximality of the Lebesgue measure. Let $\psi$ be
the maximal irreducibility measure of the chain $\{X_n\}_{n\geq0}$.
Hence, $\lambda\ll\psi.$ Let us show that $\psi\ll\lambda.$ If that
would not be the case, that is, if there would exist $B\in
\mathcal{B}(\R)$ such that $\lambda(B)=0$ and $\psi(B)>0$, then by
irreducibility of the chain $\{X_n\}_{n\geq0}$, for every $x\in\R$
there would exist $n\in\N$ such that
\[
p^{n}(x,B)=\int_{\R}p(x,\mathrm{d}x_1)\int
_{\R}p(x_1,\mathrm{d}x_2)\cdots\int
_{\R
}p(x_{n-2},\mathrm{d}x_{n-1})\int
_{B-x_{n-1}}f_{x_{n-1}}(x_n)\,\mathrm{d}x_n>0.
\]
But, since $\int_{B-x}f_{x}(y)\,\mathrm{d}y=0$, for every $x\in\R$, because
$\lambda(B)=0$, we have $p^{n}(x,B)=0$.
\end{pf}

\begin{definition}\label{d22}Let $\{Y_n\}_{n\geq0}$ be a Markov chain
on $(\R, \mathcal{B}(\R)).$
\begin{enumerate}[(iii)]
\item[(i)] A set $C\in B(\R)$ is called a \emph{$\nu_n$-small set} if
there exist $n\in\N$
and a nontrivial measure $\nu_n$ on $\mathcal{B}(\R)$ such that
for every $B\in\mathcal{B}(\R)$ and for every $x\in C$ we have
%
\begin{equation}
\label{eq21}p^{n}(x,B)\geq\nu_n(B).
\end{equation}

\item[(ii)] The $\psi$-irreducible Markov chain $\{Y_n\}_{n\geq0}$ is
called \emph{aperiodic}
if for some small set $C$ with $\psi(C)>0$, $1$ is the greatest
common divisor of all values $m\in\N$ for which (\ref{eq21}) holds
for $\nu_m=\delta_m\nu_n$, where $n\in\N$ is such that $C$ is
$\nu_n$-small set with $\nu_n(C)>0$ and $\delta_m>0$.
\item[(iii)] Let $C\in\mathcal{B}(\R)$. If there exist a probability
measure $a=\{a(n)\}_{n\geq0}$
on $\ZZ_+$ and a nontrivial measure $\nu_a$ on $\mathcal{B}(\R)$
such that
\[
\sum_{n=0}^{\infty}a(n)p^{n}(x,B)
\geq\nu_a(B)
\]
holds for every $x\in C$ and every $B\in
\mathcal{B}(\R)$, then the set $C$ is called $\nu_a$\emph{-petite
set.}
\end{enumerate}
\end{definition}

\begin{proposition}\label{p23}Conditions \emph{(C1)--(C4)} imply that
for the chain $\{X_n\}_{n\geq0}$ every bounded Borel set
$C\subseteq[-k,k]^{c}$ is a $\nu_2$-small set for some nontrivial
measure $\nu_2$.
\end{proposition}
\begin{pf} By (C3), there exists $k>0$, such that for all
$0<\varepsilon<1$
there exists $y_{\varepsilon}\geq k\vee1$, such that for all
$|y|\geq y_{\varepsilon}$ we have
\[
\biggl\llvert f_x(y)\frac{|y|^{\alpha(x)+1}}{c(x)}-1\biggr\rrvert <\varepsilon
\]
for all $x\in[-k,k]^{c}$. Let $C\subseteq(-\infty,-k]$ be a bounded
Borel set. Let $x\in C$ and $B\in\mathcal{B}(\R)$ be arbitrary.
Similarly as in Proposition \ref{p21}, we have
\begin{eqnarray*}
p^{2}(x,B)&=&\int_{\R}f_x(y-x)\int
_{B-y}f_y(z)\,\mathrm{d}z\,\mathrm{d}y\geq \int_{y_{\varepsilon}}^{2y_{\varepsilon}}f_x(y-x)
\int_{(B-y)\cap
(-\infty
,-y_{\varepsilon})}f_y(z)\,\mathrm{d}z\,\mathrm{d}y
\\
&>&(1-\varepsilon)^{2} \Bigl(\inf_{x\in
C}c(x) \Bigr) \Bigl(
\inf_{y_{\varepsilon}\leq y\leq2
y_{\varepsilon}}c(y) \Bigr)\int_{y_{\varepsilon}}^{2y_{\varepsilon
}}(y-a)^{-3}
\int_{(B-y)\cap(-\infty,-y_{\varepsilon})}|z|^{-3}\,\mathrm{d}z\,\mathrm{d}y,
\end{eqnarray*}
where $a:=\inf C.$ Now, by condition (C4), the measure
\[
\nu_2(B):=(1-\varepsilon)^{2} \Bigl(\inf_{x\in C}c(x)
\Bigr) \Bigl(\inf_{y_{\varepsilon}\leq y\leq2
y_{\varepsilon}}c(y) \Bigr)\int_{y_{\varepsilon}}^{2y_{\varepsilon
}}(y-a)^{-3}
\int_{(B-y)\cap(-\infty,-y_{\varepsilon})}|z|^{-3}\,\mathrm{d}z\,\mathrm{d}y
\]
is a nontrivial measure. Therefore, the set $C$ is a $\nu_2$-small
set. Similarly, we deduce that a bounded Borel set
$C\subseteq[k,\infty)$ is a $\nu_2$-small for some nontrivial
measure $\nu_2$.
\end{pf}

\begin{proposition}\label{p24}
Under conditions \emph{(C1)--(C4)}, the chain $\{X_n\}_{n\geq0}$ is an
aperiodic chain.
\end{proposition}
\begin{pf}
From the previous proposition, we know that every bounded Borel set
$C\subseteq[-k,k]^{c}$ is a $\nu_2$-small set. Let us show that
there exists a $\nu_2$-small set $C\subseteq[-k,k]^{c}$ which is
also a $\nu_3$-small set with $\nu_3=\delta_3\nu_2$, for some
$\delta_3>0$.
Let $C=[-4y_{\varepsilon}-k, -k]$, where $\varepsilon$ and
$y_{\varepsilon}$ are given as in the previous proposition. The set
$C$ is a $\nu_2$-small set. Let us show that
\[
\inf_{x\in C}p(x,C)>0.
\]
Then, by \cite{meyn}, Proposition 5.2.4, $C$ is a $\nu_3$-small set,
where $\nu_3$ is a multiple of $\nu_2$. Similarly as in Proposition
\ref
{p21}, we have
\begin{eqnarray*}
p(x,C)&=&\int_{C-x}f_x(y)\,\mathrm{d}y\geq\int
_{(C-x)\cap(-\infty
,-y_{\varepsilon})\cup(C-x)\cap(
y_{\varepsilon},\infty)}f_x(y)\,\mathrm{d}y
\\
&>&(1-\varepsilon) \Bigl(\inf_{x\in C}c(x) \Bigr)\inf_{x\in
C}\int
_{(C-x)\cap(-\infty,-y_{\varepsilon})\cup(C-x)\cap(
y_{\varepsilon},\infty)}|y|^{-3}\,\mathrm{d}y>0.
\end{eqnarray*}
\upqed\end{pf}

The following result is a
consequence of \cite{meyn}, Proposition 5.5.2 and Theorem 5.5.7.
%
\begin{proposition}\label{p25} Conditions \emph{(C1)--(C4)} imply that for
the chain $\{X_n\}_{n\geq0}$,
a Borel set is a small set if and only if it is a petite set.
\end{proposition}

Since conditions (C3), (C4) and (C5) consider compact sets, we
get the following result which is essential in proving
Theorems \ref{tm13} and \ref{tm14}.

\begin{proposition}\label{p26}Conditions \emph{(C1)--(C5)} imply that for the
chain $\{X_n\}_{n\geq0}$, every bounded Borel set is a small set.
\end{proposition}
\begin{pf}
From Proposition \ref{p23}, we know that every bounded Borel set
$C\subseteq[-k,k]^{c}$ is a small set. By \cite{meyn}, Proposition
5.5.5, it is enough to show that $[-k,k]$ is a small set. Let
$C\subseteq(-\infty,-k]$ be a bounded Borel set, that is, a small
set. Let $0<\varepsilon<1$ be arbitrary and let $y_{\varepsilon}\geq
(k\vee l\vee1)$ (recall that $l$ is defined in condition (C5)) be
such that for all $|y|\geq y_{\varepsilon}$ we have
\[
\biggl\llvert f_x(y)\frac{|y|^{\alpha(x)+1}}{c(x)}-1\biggr\rrvert <\varepsilon
\]
for all $x\in[-k,k]^{c}$. Then, similarly as in Proposition
\ref{p21}, for every $x\in[-k,k]$, we have
\begin{eqnarray*}
p^{2}(x,C)&=&\int_{\R}f_x(y-x)\int
_{C-y}f_y(z)\,\mathrm{d}z\,\mathrm{d}y\geq \int_{y_{\varepsilon}}^{2y_{\varepsilon}}f_x(y-x)
\int_{(C-y)\cap
(-\infty
,-y_{\varepsilon})}f_y(z)\,\mathrm{d}z\,\mathrm{d}y
\\
&>&(1-\varepsilon) \Bigl(\inf_{y_{\varepsilon}\leq y\leq2
y_{\varepsilon}}c(y) \Bigr)\inf_{x\in[-k,k]}
\biggl(\int_{[y_{\varepsilon
},2y_{\varepsilon}]-x}f_x(y)\,\mathrm{d}y \biggr) \biggl(\int
_{C-2y_{\varepsilon
}-k}|z|^{-3}\,\mathrm{d}z \biggr).
\end{eqnarray*}
Now, using condition (C5), we have that $p^{2}(x,C)>0$. Therefore,
by \cite{meyn}, Proposition 5.2.4, the set $[-k,k]$ is a small
set, that is, every bounded Borel set is a small set.
\end{pf}

\section{\texorpdfstring{Proof of Theorem \protect\ref{tm13}}{Proof of Theorem 1.3}}\label{sec3}

In this section, we give a proof of Theorem \ref{tm13}. Before the
proof we recall several special functions we need. The Digamma
function is a function defined by
$\Psi(z):=\frac{\Gamma'(z)}{\Gamma(z)},$ for $z\in\CC,$ $\Re(z)>0$, where $\Gamma(z)$ is the Gamma function.
%
\begin{lemma}\label{l31}Let $a>0$ be an arbitrary real number. Then
\[
\int_{1}^{\infty}\frac{\mathrm{d}y}{y^{a}(1+y)}=
\frac{1}{2} \biggl(\Psi \biggl(\frac{a+1}{2} \biggr)-\Psi \biggl(
\frac{a}{2} \biggr) \biggr).
\]
\end{lemma}
\begin{pf}
From \cite{abramowitz}, formula 6.3.22, we have
\[
\Psi(z)=\int_{0}^{1}\frac{1-x^{z-1}}{1-x}\,\mathrm{d}x-
\gamma,
\]
for $\Re (z)>0$, where
$\gamma$ is Euler's constant. Then
\[
\Psi \biggl(\frac{a+1}{2} \biggr)-\Psi \biggl(\frac{a}{2} \biggr)=
\int_{0}^{1}\frac{x^{{a}/{2}-1}-x^{({a+1})/{2}-1}}{1-x}\,\mathrm{d}x.
\]
The claim follows by change of variables $x=y^{-2}$.
\end{pf}
The Gauss hypergeometric function is defined by the formula
%
\begin{equation}
\label{eq41}{}_2F_1(a,b,c;z):=\sum
_{n=0}^{\infty}\frac
{(a)_n(b)_n}{(c)_n}\frac{z^{n}}{n!}
\end{equation}
for $a,b,c,z\in\CC$, $c\notin\ZZ_-$, where for $w\in\CC$ and
$n\in\ZZ_+$, $(w)_n$ is defined by
\[
(w)_0=1\quad\mbox{and}\quad(w)_n=w(w+1)\cdots(w+n-1).
\]
The series (\ref{eq41}) absolutely converges on $|z|<1$,
absolutely converges on $|z|\leq1$ when $\Re (c-a-b)>0$, conditionally converges on $|z|\leq1$, except for $z=1$, when
$-1<\Re (c-b-a)\leq0$ and diverges when
$\Re (c-b-a)\leq -1$. In the case when $\Re
(c)>\Re (b)>0$, it can be analytically continued
on $\CC\setminus(1,\infty)$ by the formula
%
\begin{equation}
\label{eq42}{}_2F_1(a,b,c;z)=\frac{\Gamma(c)}{\Gamma
(b)\Gamma
(c-b)}\int
_{0}^{1}t^{b-1}(1-t)^{c-b-1}(1-tz)^{-a}\,\mathrm{d}t.
\end{equation}

The incomplete Beta function is defined by the formula
%
\begin{equation}
\label{eq43} B(x;z,w):=\int_{0}^{x}t^{z-1}(1-t)^{w-1}\,\mathrm{d}t
\end{equation}
for
$x\in[0,1]$, $\Re (z)>0$ and $\Re
(w)>0$. When $x=1$, the function $B(1;z,w)$ is called the Beta
function and
%
\begin{equation}
\label{eq44} B(1;z,w)=\frac{\Gamma(z)\Gamma(w)}{\Gamma(z+w)}.
\end{equation}

We need the following technical lemma.
%
\begin{lemma}\label{l32} Let $\alpha\dvtx \R\longrightarrow(1,2)$ be an
arbitrary function. Then for every $R\geq0$ we have
\[
\lim_{|x|\longrightarrow\infty}\frac{1}{2-\alpha(x)} \biggl(1- \biggl(\frac
{|x|}{|x|+R}
\biggr)^{2-\alpha(x)} \biggr)=0.
\]
\end{lemma}
\begin{pf} Let
$0<\varepsilon<1$ be arbitrary. Since
\[
\frac{1}{x} \bigl(1- (1-\varepsilon )^{x} \bigr)\leq-\ln (1-
\varepsilon)
\]
for all $x\in(0,1]$, we have
\[
0\leq\operatorname{\lim\sup}\limits_{|x|\longrightarrow\infty}\frac
{1}{2-\alpha
(x)} \biggl(1- \biggl(
\frac{|x|}{|x|+R} \biggr)^{2-\alpha(x)} \biggr)\leq \operatorname{\lim\sup}\limits_{|x|\longrightarrow\infty}
\frac{1- (1-\varepsilon
)^{2-\alpha(x)}}{2-\alpha(x)}\leq -\ln(1-\varepsilon).
\]
By letting
$\varepsilon\longrightarrow0$,
we have the
claim.
\end{pf}

\begin{pf*}{Proof of Theorem \ref{tm13}} The proof is divided in
four steps.

\textit{Step} 1. In the first step, we explain our strategy of the
proof. Let us define the function $V\dvtx \R\longrightarrow\R_+$ by the
formula
\[
V(x):=\ln\bigl(1+|x|\bigr).
\]
From Proposition \ref{p26}, the set $C_V(r)=\{y\dvt V(y)\leq r\}$ is
a petite set for all $r<\infty.$ We will show that there exists
$r_0>0$, big enough, such that $\int_{\R}p(x,\mathrm{d}y)V(y)-V(x)\leq0$ for
all $x\in C_V^{c}(r_0).$ Then, the desired result will follow from
\cite{meyn}, Theorem 8.4.2. Since $C_V(r)\uparrow\R$, when
$r\longrightarrow\infty$, it is enough to show that
\[
\operatorname{\lim\sup}\limits_{|x|\longrightarrow\infty}\frac{(1+|x|)^{\alpha
(x)}}{c(x)} \biggl(\int_{\R}p(x,\mathrm{d}y)V(y)-V(x)
\biggr)<0.
\]
We have
%
\begin{eqnarray}
\label{eq31} \int_{\R}p(x,\mathrm{d}y)V(y)&=&\int_{\R
}f_x(y-x)V(y)\,\mathrm{d}y=
\int_{\R}f_x(y)V(y+x)\,\mathrm{d}y
\nonumber
\\[-8pt]\\[-8pt]
&=&\int_{-x}^{\infty}\ln (1+x+y)f_x(y)\,\mathrm{d}y+
\int_{-\infty}^{-x}\ln(1-x-y)f_x(y)\,\mathrm{d}y.\nonumber
\end{eqnarray}

\textit{Step} 2. In the second step, we find an appropriate upper
bound for the first summand in (\ref{eq31}). For any $x>0$, we
have
\[
\int_{-x}^{\infty}\ln(1+x+y)f_x(y)\,\mathrm{d}y =
\ln(1+x)\int_{-x}^{\infty}f_x(y)\,\mathrm{d}y+\int
_{-x}^{\infty}\ln \biggl(1+\frac{y}{1+x}
\biggr)f_x(y)\,\mathrm{d}y.
\]
Let $0<\delta<1$ be
arbitrary. By restricting $\ln(1+t)$ to intervals $(-1,-\delta),$
$[-\delta,\delta]$, $(\delta,1)$ and $[1,\infty)$, and using the
Taylor expansion of the function $\ln(1+t)$, that is,
\[
\ln(1+t)=\sum_{i=1}^{\infty}
\frac{(-1)^{i+1}}{i}t^{i}
\]
for $t\in(-1,1]$, we get
\begin{eqnarray*}
\int_{-x}^{\infty}\ln(1+x+y)f_x(y)\,\mathrm{d}y&\leq&
\ln (1+x)\int_{-x}^{\infty}f_x(y)\,\mathrm{d}y
\\
&&{}-\sum_{i=1}^{\infty}\frac{1}{i(1+x)^{i}}\int
_{\{-1-x<y<-\delta
(1+x)\}\cap
\{y+x>0\}}|y|^{i}f_x(y)\,\mathrm{d}y
\\
&&{}+\int_{\{-\delta(1+x)\leq y\leq\delta(1+x)\}\cap
\{y+x>0\}}\ln \biggl(1+\frac{y}{1+x}
\biggr)f_x(y)\,\mathrm{d}y
\\
&&{}+\sum_{i=1}^{\infty}\frac{(-1)^{i+1}}{i(1+x)^{i}}\int
_{\{\delta
(1+x)<y<1+x\}\cap
\{y+x>0\}}y^{i}f_x(y)\,\mathrm{d}y
\\
&&{}+\int_{\{y\geq1+x\}\cap
\{y+x>0\}}\ln \biggl(1+\frac{y}{1+x}
\biggr)f_x(y)\,\mathrm{d}y.
\end{eqnarray*}
Furthermore, by taking $x>\frac{\delta}{1-\delta}$ we get
\begin{eqnarray*}
\int_{-x}^{\infty}\ln(1+x+y)f_x(y)\,\mathrm{d}y&\leq&
\ln (1+x)\int_{-x}^{\infty}f_x(y)\,\mathrm{d}y-\sum
_{i=1}^{\infty}\frac{1}{i(1+x)^{i}}\int
_{-x}^{-\delta
(1+x)}|y|^{i}f_x(y)\,\mathrm{d}y
\\
&&{}+\int_{-\delta(1+x)}^{\delta(1+x)}\ln \biggl(1+\frac
{y}{1+x}
\biggr)f_x(y)\,\mathrm{d}y
\\
&&{}+\sum_{i=1}^{\infty}\frac
{(-1)^{i+1}}{i(1+x)^{i}}\int
_{\delta
(1+x)}^{1+x}y^{i}f_x(y)\,\mathrm{d}y
\\
&&{}+\int_{1+x}^{\infty}\ln \biggl(1+\frac{y}{1+x}
\biggr)f_x(y)\,\mathrm{d}y.
\end{eqnarray*}
Let us put
\begin{eqnarray*}
U^{\delta}_1(x)&:=&-\frac{1}{1+x}\int
_{\delta
(1+x)}^{x}yf_x(-y)\,\mathrm{d}y+
\frac{1}{1+x}\int_{\delta
(1+x)}^{1+x}yf_x(y)\,\mathrm{d}y,
\\
U^{\delta}_2(x)&:=&-\frac{1}{2(1+x)^{2}}\int
_{\delta
(1+x)}^{x}y^{2}f_x(-y)\,\mathrm{d}y-
\frac{1}{2(1+x)^{2}}\int_{\delta
(1+x)}^{1+x}y^{2}f_x(y)\,\mathrm{d}y,
\\
U^{\delta}_3(x)&:=&-\sum_{i=3}^{\infty}
\frac
{1}{i(1+x)^{i}}\int_{\delta
(1+x)}^{x}y^{i}f_x(-y)\,\mathrm{d}y+
\sum_{i=3}^{\infty}\frac
{(-1)^{i+1}}{i(1+x)^{i}}\int
_{\delta
(1+x)}^{1+x}y^{i}f_x(y)\,\mathrm{d}y,
\\
U^{\delta}_4(x)&:=&\int_{-\delta(1+x)}^{\delta(1+x)}
\ln \biggl(1+\frac
{y}{1+x} \biggr)f_x(y)\,\mathrm{d}y \quad\mbox{and}
\\
U_5(x)&:=&\int_{1+x}^{\infty}\ln \biggl(1+
\frac{y}{1+x} \biggr)f_x(y)\,\mathrm{d}y
\end{eqnarray*}
for $0<\delta<1$ and $x>\frac{\delta}{1-\delta}$. Hence, we find
%
\begin{eqnarray}
\label{eq32}&&\int_{-x}^{\infty}\ln (1+x+y)f_x(y)\,\mathrm{d}y
\nonumber
\\[-8pt]\\[-8pt]
&&\quad\leq\ln(1+x)\int_{-x}^{\infty}f_x(y)\,\mathrm{d}y+
U^{\delta}_1(x)+ U^{\delta}_2(x)+U^{\delta}_3(x)+U^{\delta}_4(x)+U_5(x).\nonumber
\end{eqnarray}
Here comes the crucial step where condition (C3) is needed. In the
above terms, by (C3), we can replace all the density functions
$f_x(y)$ by the functions $c(x)|y|^{-\alpha(x)-1}$ and find a more
operable upper bound in (\ref{eq32}). Let $0<\varepsilon<1$ be
arbitrary. Then, by (C3), there exists $y_{\varepsilon}\geq1$,
such that for all $|y|\geq y_{\varepsilon}$
\[
\biggl\llvert f_x(y)\frac{|y|^{\alpha(x)+1}}{c(x)}-1\biggr\rrvert <\varepsilon
\]
for all $x\in[-k,k]^{c}.$
Let
$x> (k\vee\frac{y_{\varepsilon}-\delta}{\delta}\vee\frac
{\delta
}{1-\delta} ).$
By a straightforward calculation, we have
\begin{eqnarray*}
U^{\delta}_1(x)&<&-\frac{(1-\varepsilon)c(x)}{(\alpha
(x)-1)(1+x)^{\alpha(x)}} \biggl(
\delta^{-\alpha(x)+1}- \biggl(\frac
{x}{1+x} \biggr)^{-\alpha(x)+1} \biggr)
\\
&&{}+\frac{(1+\varepsilon
)c(x)}{(\alpha(x)-1)(1+x)^{\alpha(x)}}\frac{\delta-\delta^{\alpha
(x)}}{\delta^{\alpha(x)}}=:U_1^{\delta,\varepsilon}(x),
\\
U^{\delta}_2(x)&<&-\frac{(1-\varepsilon)c(x)}{(1+x)^{\alpha
(x)}}\frac
{1}{2(2-\alpha(x))}
\biggl( \biggl(\frac{x}{1+x} \biggr)^{2-\alpha
(x)}-\delta^{2-\alpha(x)}
\biggr)
\\
&&{}-\frac{(1-\alpha)c(x)}{(1+x)^{\alpha(x)}}\frac{1}{2(2-\alpha
(x))}\frac
{\delta^{\alpha(x)}-\delta^{2}}{\delta^{\alpha(x)}}=:U_2^{\delta
,\varepsilon}(x),
\\
U_3^{\delta}(x)&<&-\frac{(1-\varepsilon)c(x)}{(1+x)^{\alpha
(x)}}\sum
_{i=3}^{\infty}\frac{1}{i(i-\alpha(x))} \biggl( \biggl(
\frac{x}{1+x} \biggr)^{i-\alpha(x)}-\delta^{i-\alpha(x)} \biggr)
\\
&&{}+\frac{c(x)}{(1+x)^{\alpha(x)}}\sum_{i=3}^{\infty}
\biggl(\frac{(-1)^{i+1}(1+(-1)^{i+1}\varepsilon)}{i(i-\alpha(x))}\frac
{\delta^{\alpha(x)}-\delta^{i}}{\delta^{\alpha(x)}} \biggr)=:U_3^{\delta
,\varepsilon}(x)
\quad\mbox{and}
\\
U_5(x)&<&(1+\varepsilon)c(x)\int_{1+x}^{\infty}
\ln \biggl(1+\frac
{y}{1+x} \biggr)\frac{1}{y^{\alpha(x)+1}}\,\mathrm{d}y=:U^{\varepsilon}_5(x).
\end{eqnarray*}
Hence, from (\ref{eq32}), we get
%
\begin{eqnarray}
\label{eq33}&&\int_{-x}^{\infty}\ln (1+x+y)f_x(y)\,\mathrm{d}y
\nonumber
\\[-8pt]\\[-8pt]
&&\quad<\ln(1+x)\int_{-x}^{\infty}f_x(y)\,\mathrm{d}y
+U^{\delta,\varepsilon}_1(x)+ U^{\delta,\varepsilon}_2(x)+U^{\delta,\varepsilon}_3(x)+U^{\delta
}_4(x)+U_5^{\varepsilon}(x).\nonumber
\end{eqnarray}

\textit{Step} 3. In the third step, we find an appropriate upper
bound for the second summand in (\ref{eq31}). We have
\[
\int_{\infty}^{-x}\ln(1-x-y)f_x(y)\,\mathrm{d}y=\ln
(x-1)\int_{\infty}^{-x} f_x(y)\,\mathrm{d}y+\int
_{\infty}^{-x}\ln \biggl(-1-\frac{y}{x-1}
\biggr)f_x(y)\,\mathrm{d}y .
\]
Let
$x> (k\vee\frac{y_{\varepsilon}-\delta}{\delta}\vee\frac
{\delta
}{1-\delta} ).$
Then, again by (C3),
\begin{eqnarray*}
\int_{\infty}^{-x}\ln(1-x-y)f_x(y)\,\mathrm{d}y &<&
\ln(x-1)\int_{\infty}^{-x} f_x(y)\,\mathrm{d}y
\\
&&{}+c(x) (1-\varepsilon)\int_{x}^{2x-2}\ln \biggl(-1+
\frac
{y}{x-1} \biggr)\frac{\mathrm{d}y}{|y|^{\alpha(x)+1}}
\\
&&{}+c(x) (1+\varepsilon )\int_{2x-2}^{\infty}\ln \biggl(-1+
\frac{y}{x-1} \biggr)\frac
{\mathrm{d}y}{|y|^{\alpha
(x)+1}}
\\
&=&\ln(x-1)\int_{\infty}^{-x} f_x(y)\,\mathrm{d}y
\\
&&{}+c(x) (1-\varepsilon)\int_{x}^{\infty}\ln \biggl(-1+
\frac
{y}{x-1} \biggr)\frac{\mathrm{d}y}{|y|^{\alpha(x)+1}}
\\
&&{}+2\varepsilon c(x)\int_{2x-2}^{\infty}\ln \biggl(-1+
\frac{y}{x-1} \biggr)\frac
{\mathrm{d}y}{|y|^{\alpha(x)+1}}.
\end{eqnarray*}
Let us put
\begin{eqnarray*}
U_6^{\varepsilon}(x)&:=&c(x) (1-\varepsilon)\int
_{x}^{\infty
}\ln \biggl(-1+\frac{y}{x-1} \biggr)
\frac{\mathrm{d}y}{|y|^{\alpha(x)+1}}
\\
&&{}+2\varepsilon c(x)\int_{2x-2}^{\infty}\ln \biggl(-1+
\frac{y}{x-1} \biggr)\frac
{\mathrm{d}y}{|y|^{\alpha(x)+1}}.
\end{eqnarray*}
We have
%
\begin{equation}
\label{eq34}\int_{\infty}^{-x}\ln(1-x-y)f_x(y)\,\mathrm{d}y
<\ln(x-1)\int_{\infty}^{-x}f_x(y)\,\mathrm{d}y+U_6^{\varepsilon}(x).
\end{equation}

\textit{Step} 4. In the fourth step, we prove
\[
\operatorname{\lim\sup}\limits_{x\longrightarrow\infty}\frac{(1+x)^{\alpha
(x)}}{c(x)} \biggl(\int_{\R}p(x,\mathrm{d}y)V(y)-V(x)
\biggr)<0.
\]
By combining
(\ref{eq31}), (\ref{eq33}) and (\ref{eq34}), we have
\begin{eqnarray*}
\int_{\R}p(x,\mathrm{d}y)V(y)&<&U_0(x)+U^{\delta,\varepsilon}_1(x)+
U^{\delta,\varepsilon}_2(x)+U^{\delta,\varepsilon}_3(x)\\
&&{}+U^{\delta
}_4(x)+U_5^{\varepsilon}(x)+U_6^{\varepsilon}(x),
\end{eqnarray*}
where
\begin{eqnarray*}
U_0(x)&=&\ln(1+x)\int_{-x}^{\infty}f_x(y)\,\mathrm{d}y+
\ln (x-1)\int_{-\infty}^{-x}f_x(y)\,\mathrm{d}y
\\
&=&\ln(1+x)-\ln(1+x)\int_{-\infty}^{-x}f_x(y)\,\mathrm{d}y+
\ln(x-1)\int_{-\infty
}^{-x}f_x(y)\,\mathrm{d}y
\\
&<&\ln(1+x)=V(x).
\end{eqnarray*}
Hence, we have
%
\begin{equation}
\label{eq35}
\int_{\R}p(x,\mathrm{d}y)V(y)-V(x)<U^{\delta
,\varepsilon}_1(x)+
U^{\delta,\varepsilon}_2(x)+U^{\delta,\varepsilon}_3(x)+U^{\delta
}_4(x)+U_5^{\varepsilon}(x)+U_6^{\varepsilon}(x).
\end{equation}
In the rest of the fourth step, we prove
\begin{eqnarray*}
&&\operatorname{\lim\sup}\limits_{x\longrightarrow\infty}\frac{(1+x)^{\alpha
(x)}}{c(x)} \biggl(\int_{\R}p(x,\mathrm{d}y)V(y)-V(x)
\biggr)
\\
&&\quad<\operatorname{\lim\sup}\limits_{\delta
\longrightarrow0}\operatorname{\lim\sup}\limits_{\varepsilon\longrightarrow0}\operatorname{\lim\sup}\limits_{x\longrightarrow\infty}
\frac{(1+x)^{\alpha(x)}}{c(x)}U^{\delta
,\varepsilon}_1(x)+\operatorname{\lim\sup}\limits_{\delta\longrightarrow0}
\operatorname{\lim\sup}\limits_{\varepsilon
\longrightarrow0}\operatorname{\lim\sup}\limits_{x\longrightarrow\infty}\frac
{(1+x)^{\alpha
(x)}}{c(x)}U^{\delta,\varepsilon}_2(x)
\\
&&\qquad {}+\operatorname{\lim\sup}\limits_{\delta
\longrightarrow0}\operatorname{\lim\sup}\limits_{\varepsilon\longrightarrow0}\operatorname{\lim\sup}\limits_{x\longrightarrow\infty}
\frac{(1+x)^{\alpha(x)}}{c(x)}U^{\delta
,\varepsilon}_3(x)+\operatorname{\lim\sup}\limits_{\varepsilon\longrightarrow0}
\operatorname{\lim\sup}\limits_{x\longrightarrow\infty}\frac{(1+x)^{\alpha
(x)}}{c(x)}U_5^{\varepsilon
}(x)
\\
&&\qquad {}+\operatorname{\lim\sup}\limits_{\varepsilon\longrightarrow0}\operatorname{\lim\sup}\limits_{x\longrightarrow\infty}\frac{(1+x)^{\alpha
(x)}}{c(x)}U_6^{\varepsilon
}(x)+R(
\alpha)\leq0.
\end{eqnarray*}
Recall
that $\alpha=\operatorname{\lim\inf}_{|x|\longrightarrow\infty}\alpha(x)>1$,
\[
R(\alpha)=\sum_{i=1}^{\infty}
\frac{1}{i(2i-\alpha)}-\frac
{\ln
2}{\alpha}-\frac{1}{2\alpha} \biggl(\Psi \biggl(
\frac{\alpha
+1}{2} \biggr)-\Psi \biggl(\frac{\alpha}{2} \biggr) \biggr)
\]
and
\[
\operatorname{\lim\sup}\limits_{\delta\longrightarrow0}\operatorname{\lim\sup}\limits_{x\longrightarrow\infty
}\frac
{(1+x)^{\alpha(x)}}{c(x)}U_4^{\delta}(x)<R(
\alpha)
\]
when $\alpha<2$, and the above limit is finite when $\alpha=2$
(assumption (\ref{eq13})). We have
%
\begin{eqnarray}
\label{eq37}&&\operatorname{\lim\sup}\limits_{\delta\longrightarrow0}
\operatorname{\lim\sup}\limits_{\varepsilon\longrightarrow0}\operatorname{\lim\sup}\limits_{x\longrightarrow\infty
}
\frac
{(1+x)^{\alpha(x)}}{c(x)}U^{\delta,\varepsilon}_1(x)
\nonumber
\\
&&\quad=\operatorname{\lim\sup}\limits_{\delta\longrightarrow0}\operatorname{\lim\sup}\limits_{\varepsilon\longrightarrow
0}\operatorname{\lim\sup}\limits_{x\longrightarrow\infty} \biggl[-
\frac{1-\varepsilon}{\alpha
(x)-1} \biggl(\delta^{-\alpha(x)+1}- \biggl(\frac{x}{1+x}
\biggr)^{-\alpha
(x)+1} \biggr)
\nonumber
\\
&&\hspace*{113pt}{} +\frac{1+\varepsilon}{\alpha(x)-1}\frac{\delta
-\delta^{\alpha(x)}}{\delta^{\alpha(x)}} \biggr]
\nonumber
\\
&&\quad=\operatorname{\lim\sup}\limits_{\delta\longrightarrow0}\operatorname{\lim\sup}\limits_{\varepsilon
\longrightarrow
0}\operatorname{\lim\sup}\limits_{x\longrightarrow\infty} \biggl[-
\frac{1-\varepsilon
}{\alpha
(x)-1} \biggl(\frac{\delta-\delta^{\alpha(x)}}{\delta^{\alpha
(x)}}+1- \biggl(\frac{x}{1+x}
\biggr)^{-\alpha(x)+1} \biggr)
\\
&&\hspace*{113pt}{} +\frac
{1+\varepsilon}{\alpha(x)-1}\frac{\delta-\delta^{\alpha
(x)}}{\delta^{\alpha(x)}} \biggr]
\nonumber
\\
&&\quad=\operatorname{\lim\sup}\limits_{\delta\longrightarrow0}\operatorname{\lim\sup}\limits_{\varepsilon
\longrightarrow
0}\operatorname{\lim\sup}\limits_{x\longrightarrow\infty} \biggl[
\frac{2\varepsilon
}{\alpha
(x)-1}\frac{\delta-\delta^{\alpha(x)}}{\delta^{\alpha(x)}}-\frac
{1-\varepsilon}{\alpha(x)-1} \biggl(1- \biggl(
\frac{x}{1+x} \biggr)^{-\alpha
(x)+1} \biggr) \biggr]
\nonumber
\\
&&\quad=\operatorname{\lim\sup}\limits_{x\longrightarrow\infty} \biggl[\frac{1}{\alpha
(x)-1} \biggl( \biggl(
\frac{x}{x+1} \biggr)^{-\alpha(x)+1}-1 \biggr) \biggr]=0.\nonumber
\end{eqnarray}
In the last two equalities, we use the assumption
$\operatorname{\lim\inf}_{|x|\longrightarrow\infty}\alpha(x)>1.$ From
Lemma \ref{l32}, we have
%
\begin{eqnarray}
\label{eq38}&&\operatorname{\lim\sup}\limits_{\delta\longrightarrow0}\operatorname{\lim\sup}\limits_{\varepsilon\longrightarrow0}\operatorname{\lim\sup}\limits_{x\longrightarrow\infty
}
\frac
{(1+x)^{\alpha(x)}}{c(x)}U^{\delta,\varepsilon}_2(x)
\nonumber
\\
&&\quad=\operatorname{\lim\sup}\limits_{\delta\longrightarrow0}\operatorname{\lim\sup}\limits_{\varepsilon\longrightarrow
0}\operatorname{\lim\sup}\limits_{x\longrightarrow\infty} \biggl[-
\frac{1-\varepsilon}{2(2-\alpha
(x))} \biggl( \biggl(\frac{x}{1+x} \biggr)^{2-\alpha(x)}-
\delta^{2-\alpha
(x)} \biggr)
\nonumber
\\
&&\hspace*{113pt} {}-\frac{1-\varepsilon}{2(2-\alpha
(x))}\frac{\delta^{\alpha(x)}-\delta^{2}}{\delta^{\alpha
(x)}} \biggr]
\nonumber
\\
&&\quad=\operatorname{\lim\sup}\limits_{\delta\longrightarrow0}\operatorname{\lim\sup}\limits_{x\longrightarrow\infty
} \biggl[-\frac{1}{2(2-\alpha(x))}
\biggl( \biggl(\frac{x}{1+x} \biggr)^{2-\alpha
(x)}+\frac{\delta^{\alpha(x)}-\delta^{2}}{\delta^{\alpha
(x)}}-1
\biggr)
\\
&&\hspace*{82pt} {}-\frac{1}{2(2-\alpha(x))}\frac{\delta^{\alpha
(x)}-\delta^{2}}{\delta^{\alpha(x)}} \biggr]
\nonumber
\\
&&\quad=\operatorname{\lim\sup}\limits_{\delta\longrightarrow0}\operatorname{\lim\sup}\limits_{x\longrightarrow\infty
} \biggl[-\frac{1}{2-\alpha(x)}
\frac{\delta^{\alpha(x)}-\delta^{2}}{\delta^{\alpha(x)}} \biggr]\nonumber\\
&&\quad\leq\cases{ %
\displaystyle -\frac{1}{2-\alpha}, &\quad $\alpha<2$,
\cr
-\infty,&\quad $\alpha=2.$}\nonumber
\end{eqnarray}
Using the dominated convergence theorem, we have
%
\begin{eqnarray}\label{eq39}
&&\hspace*{-12pt}\operatorname{\lim\sup}\limits_{\delta\longrightarrow
0}\operatorname{\lim\sup}\limits_{\varepsilon\longrightarrow0}\operatorname{\lim\sup}\limits_{x\longrightarrow\infty
}
\frac
{(1+x)^{\alpha(x)}}{c(x)}U^{\delta,\varepsilon}_3(x)
\nonumber
\\
&&\hspace*{-12pt}\quad=\operatorname{\lim\sup}\limits_{\delta\longrightarrow0}\operatorname{\lim\sup}\limits_{\varepsilon
\longrightarrow0}
\operatorname{\lim\sup}\limits_{x\longrightarrow\infty} \Biggl[-(1-\varepsilon )\sum_{i=3}^{\infty}
\frac{1}{i(i-\alpha(x))} \biggl( \biggl(\frac{x}{1+x} \biggr)^{i-\alpha(x)}-
\delta^{i-\alpha(x)} \biggr)
\nonumber\\
&&\hspace*{-12pt}\hspace*{114pt}{}+\sum_{i=3}^{\infty}\frac{(\varepsilon
+(-1)^{i+1})}{i(i-\alpha(x))}
\frac{\delta^{\alpha(x)}-\delta^{i}}{\delta^{\alpha(x)}} \Biggr] %
\nonumber
\\
&&\hspace*{-12pt}\quad=\operatorname{\lim\sup}\limits_{\delta\longrightarrow0}\operatorname{\lim\sup}\limits_{\varepsilon
\longrightarrow0}
\operatorname{\lim\sup}\limits_{x\longrightarrow\infty} \Biggl[ \sum_{i=3}^{\infty}
\frac{- ({x}/({1+x}) )^{i-\alpha
(x)}+\delta^{i-\alpha(x)}+(-1)^{i+1}-(-1)^{i+1}\delta^{i-\alpha
(x)}}{i(i-\alpha(x))}\nonumber
\\[-8pt]\\[-8pt]
&&\hspace*{-12pt}\hspace*{114pt}{}+\varepsilon\sum_{i=3}^{\infty}
\frac{ (
{x}/({1+x}) )^{i-\alpha(x)}-\delta^{i-\alpha(x)}+1-\delta^{i-\alpha
(x)}}{i(i-\alpha(x))} \Biggr]\qquad %
\nonumber
\\
&&\hspace*{-12pt}\quad=\operatorname{\lim\sup}\limits_{\delta\longrightarrow0}\operatorname{\lim\sup}\limits_{x\longrightarrow\infty
}\sum
_{i=3}^{\infty}\frac{- ({x}/({1+x})
)^{i-\alpha(x)}+\delta^{i-\alpha(x)}+(-1)^{i+1}-(-1)^{i+1}\delta^{i-\alpha(x)}}{i(i-\alpha(x))}
\nonumber
\\
&&\hspace*{-12pt}\quad=\operatorname{\lim\sup}\limits_{\delta\longrightarrow0}\operatorname{\lim\sup}\limits_{x\longrightarrow\infty
}\sum
_{i=3}^{\infty} \biggl(\frac{- ({x}/({1+x})
)^{i-\alpha(x)}+(-1)^{i+1}}{i(i-\alpha(x))}+
\frac{\delta^{i-\alpha
(x)}-(-1)^{i+1}\delta^{i-\alpha(x)}}{i(i-\alpha(x))} \biggr)
\nonumber
\\
&&\hspace*{-12pt}\quad\leq-\sum_{i=2}^{\infty}
\frac{2}{2i(2i-\alpha
)}=-\sum_{i=2}^{\infty}
\frac{1}{i(2i-\alpha)}.\nonumber
\end{eqnarray}
Therefore, by combining (\ref{eq37}), (\ref{eq38}) and
(\ref{eq39}) we get
%
\begin{eqnarray}
\label{eq310}
&&\operatorname{\lim\sup}\limits_{\delta\longrightarrow0}\operatorname{\lim\sup}
\limits_{\varepsilon\longrightarrow0}\operatorname{\lim\sup}\limits_{x\longrightarrow\infty}
\frac
{(1+x)^{\alpha(x)}}{c(x)} \bigl(U^{\delta,\varepsilon
}_1(x)+U^{\delta
,\varepsilon}_2(x)+U^{\delta,\varepsilon}_3(x)
\bigr)
\nonumber
\\[-8pt]\\[-8pt]
&&\quad\leq \cases{ %
\displaystyle -\sum
_{i=1}^{\infty}\frac{1}{i(2i-\alpha)}, & \quad $\alpha<2$,
\cr
-\infty,& \quad $\alpha=2$. }\nonumber
\end{eqnarray}

Now, let us calculate
\[
\operatorname{\lim\sup}\limits_{\varepsilon\longrightarrow0}\operatorname{\lim\sup}\limits_{x\longrightarrow
\infty
}\frac{(1+x)^{\alpha(x)}}{c(x)}U^{\varepsilon}_5(x).
\]
Using integration by parts formula, we get
\begin{eqnarray*}
&&\operatorname{\lim\sup}\limits_{\varepsilon\longrightarrow0}\operatorname{\lim\sup}\limits_{x\longrightarrow\infty}\frac{(1+x)^{\alpha
(x)}}{c(x)}U^{\varepsilon
}_5(x)\\&&\quad=
\operatorname{\lim\sup}\limits_{x\longrightarrow\infty}(1+x)^{\alpha(x)}\int_{1+x}^{\infty}
\ln \biggl(1+\frac{y}{1+x} \biggr)\frac{1}{y^{\alpha
(x)+1}}\,\mathrm{d}y
\\
&&\quad=\operatorname{\lim\sup}\limits_{x\longrightarrow\infty} \biggl(\frac{\ln2}{\alpha
(x)}+\frac
{1}{\alpha(x)}\int
_{1}^{\infty}\frac{\mathrm{d}y}{y^{\alpha(x)}(1+y)} \biggr).
\end{eqnarray*}
Furthermore, from Lemma \ref{l31} and the fact that the function
\[
x\longmapsto\Psi \biggl(\frac{x+1}{2} \biggr)-\Psi \biggl(
\frac
{x}{2} \biggr)
\]
is decreasing on $(0,\infty)$ (Lemma \ref{l31}) we have
%
\begin{eqnarray}\label{eq311}
&&\operatorname{\lim\sup}\limits_{\varepsilon\longrightarrow
0}\operatorname{\lim\sup}\limits_{x\longrightarrow\infty}\frac{(1+x)^{\alpha
(x)}}{c(x)}U^{\varepsilon}_5(x)
\nonumber
\\
&&\quad= \operatorname{\lim\sup}\limits_{x\longrightarrow\infty} \biggl(\frac{\ln
2}{\alpha(x)}+\frac{1}{2\alpha(x)} \biggl(
\Psi \biggl(\frac{\alpha
(x)+1}{2} \biggr)-\Psi \biggl(\frac{\alpha(x)}{2} \biggr)
\biggr) \biggr)
\\
&&\quad\leq\frac{\ln
2}{\alpha}+\frac{1}{2\alpha} \biggl(\Psi \biggl(
\frac{\alpha
+1}{2} \biggr)-\Psi \biggl(\frac{\alpha}{2} \biggr) \biggr).\nonumber
\end{eqnarray}

At the end, using integration by parts formula, we have
\begin{eqnarray*}
&&\operatorname{\lim\sup}\limits_{\varepsilon\longrightarrow0}\operatorname{\lim\sup}\limits_{x\longrightarrow\infty}\frac{(1+x)^{\alpha
(x)}}{c(x)}U_6^{\varepsilon
}(x)
\\
&&\quad=\operatorname{\lim\sup}\limits_{\varepsilon\longrightarrow0}\operatorname{\lim\sup}\limits_{x\longrightarrow
\infty}(1+x)^{\alpha(x)} \biggl[(1-
\varepsilon)\int_{x}^{\infty}\ln \biggl(-1+
\frac{y}{x-1} \biggr)\frac{1}{|y|^{\alpha(x)+1}}\,\mathrm{d}y
\\
&&\hspace*{129pt}{} +2\varepsilon \int_{2x-2}^{\infty}\ln \biggl(-1+
\frac{y}{x-1} \biggr)\frac
{1}{|y|^{\alpha
(x)+1}}\,\mathrm{d}y \biggr]
\\
&&\quad=\operatorname{\lim\sup}\limits_{\varepsilon\longrightarrow0}\operatorname{\lim\sup}\limits_{x\longrightarrow
\infty
}(1+x)^{\alpha(x)} \biggl[
\frac{1-\varepsilon}{\alpha(x)} \biggl(\frac
{1}{x^{\alpha(x)}}\ln \biggl(-1+\frac{x}{x-1}
\biggr)+\int_{x}^{\infty
}\frac{\mathrm{d}y}{y^{\alpha(x)}(y-x+1)} \biggr)
\\
&&\hspace*{129pt}{} +\frac
{2\varepsilon}{\alpha(x)}\int_{2x-2}^{\infty}
\frac{\mathrm{d}y}{y^{\alpha
(x)}(y-x+1)} \biggr]
\\
&&\quad=\operatorname{\lim\sup}\limits_{\varepsilon\longrightarrow0}\operatorname{\lim\sup}\limits_{x\longrightarrow
\infty
} \biggl[\frac{1-\varepsilon}{\alpha(x)}
\biggl(\frac{(1+x)^{\alpha
(x)}}{x^{\alpha(x)}}\ln \biggl(\frac{1}{x-1} \biggr)+
\frac
{(1+x)^{\alpha
(x)}}{(x-1)^{\alpha(x)}}\int_0^{({x-1})/{x}}
\frac{y^{\alpha
(x)-1}}{1-y}\,\mathrm{d}y \biggr)
\\
&&\hspace*{84pt}{} +\frac{2\varepsilon
}{\alpha^{2}(x)}\frac{(1+x)^{\alpha(x)}}{(x-1)^{\alpha(x)}}{}_2F_1\bigl(
\alpha(x),\alpha(x),\alpha(x)+1;-1\bigr) \biggr],
\end{eqnarray*}
where in the last equality we use (\ref{eq42}). From
(\ref{eq42}), we get
\[
{}_2F_1\bigl(\alpha(x),\alpha(x),\alpha(x)+1;-1\bigr)
\leq2\int_{0}^{1}(1+t)^{-1}\,\mathrm{d}t=\ln4,
\]
and
\[
\int_0^{({x-1})/{x}}\frac{y^{\alpha(x)-1}}{1-y}\,\mathrm{d}y\leq\int
_0^{
({x-1})/{x}}\frac{\mathrm{d}y}{1-y}=\ln x.
\]
Hence,
%
\begin{eqnarray}
\label{eq312}
&&\operatorname{\lim\sup}\limits_{\varepsilon\longrightarrow
0}\operatorname{\lim\sup}\limits_{x\longrightarrow\infty}
\frac{(1+x)^{\alpha
(x)}}{c(x)}U_6^{\varepsilon}(x)
\nonumber
\\
&&\quad\leq\operatorname{\lim\sup}\limits_{x\longrightarrow\infty}\frac{1}{\alpha(x)} \biggl(\frac
{(1+x)^{\alpha(x)}}{x^{\alpha(x)}}\ln
\biggl(-1+\frac{x}{x-1} \biggr)+\frac
{(x+1)^{\alpha(x)}}{(x-1)^{\alpha(x)}}\ln x \biggr)
\\
&&\quad\leq\operatorname{\lim\sup}\limits_{x\longrightarrow\infty}\frac{1}{\alpha(x)} \biggl(\ln \biggl(-1+
\frac{x}{x-1} \biggr)+ \biggl(1+\frac{2}{x-1} \biggr)^{2}\ln
x \biggr)=0.\nonumber
\end{eqnarray}
By combining (\ref{eq35}),
(\ref{eq310}), (\ref{eq311}) and (\ref{eq312}), we have
\[
\operatorname{\lim\sup}\limits_{x\longrightarrow\infty}\frac{(1+x)^{\alpha
(x)}}{c(x)} \biggl(\int_{\R}p(x,\mathrm{d}y)V(y)-V(x)
\biggr)<0.
\]
The case when $x<0$ is treated in the same way. Therefore,
we have proved the desired result.
\end{pf*}

\section{\texorpdfstring{Proof of Theorem \protect\ref{tm14}}{Proof of Theorem 1.4}}\label{sec4}

Let us first list some properties of the Gauss hypergeometric
function which will be needed in the proof of Theorem
\ref{tm14}:
\begin{enumerate}[(iii)]
\item[(i)] for $a,b,c,z\in\CC, c\notin\ZZ_-,$
%
\begin{equation}
_2\label{eq45}F_1(0,b,c;z) ={}_2F_1(a,0,c;z)=1;
\end{equation}
\item[(ii)] for $\Re (c-a-b)>0$, $c\notin\ZZ_-,$
%
\begin{equation}
\label{eq46}{}_2F_1(a,b,c;1)=\frac{\Gamma(c)\Gamma(c-a-b)}{\Gamma
(c-a)\Gamma(c-b)};
\end{equation}
\item[(iii)] for
$z\in\CC\setminus(1,\infty)$
%
\begin{equation}
\label{eq47}{}_2F_1(a,b,c;z)=(1-z)^{c-b-a}{}_2F_1(c-a,c-b,c;z);
\end{equation}
\item[(iv)] for $z\in\CC\setminus(0,\infty)$
%
\begin{eqnarray}\label{eq48}
{}_2F_1(a,b,c;z)&=&\frac{\Gamma(c)\Gamma(b-a)}{\Gamma(b)\Gamma
(c-a)}(-z)^{-a}{}_2F_1
\biggl(a,1-c+a,1-b+a,\frac{1}{z} \biggr)
\nonumber
\\[-8pt]\\[-8pt]
&&{}+\frac{\Gamma(c)\Gamma(a-b)}{\Gamma(a)\Gamma(c-b)}(-z)^{-b}{}_2F_1
\biggl(b,1-c+b,1-a+b,\frac{1}{z} \biggr).\nonumber
\end{eqnarray}
\end{enumerate}
For further
properties of the hypergeometric functions, the incomplete Beta
functions and the Beta function (see \cite{abramowitz}, Chapters 6 and
15).

\begin{pf*}{Proof of Theorem \ref{tm14}} The proof is divided in
three steps.

\textit{Step} 1. In the first step, we explain our strategy of the
proof. Let us define the function $V\dvtx \R\longrightarrow\R_+$ by the
formula
\[
V(x):=1-(1+|x|)^{-\beta},
\]
where $0<\beta<1-\alpha$ is arbitrary (recall that $\alpha=\operatorname{\lim\sup}_{|x|\longrightarrow\infty}\alpha(x)<1$). It is clear that
$C_V(r)\in
\mathcal{B}^{+}(\R)$ and $C^{c}_V(r)\in\mathcal{B}^{+}(\R)$, for
every $0<r<1.$
By \cite{meyn}, Theorem 8.4.3, we have to show that there exists
$0<r_0<1$ such that $\Delta V(x)\geq0$, for every $x\in
C_V^{c}(r_0).$ Since $C_V(r)\uparrow\R$, when $r\uparrow1$, it is
enough to show that
\[
\operatorname{\lim\inf}\limits_{|x|\longrightarrow\infty}\frac{\alpha
(x)|x|^{\alpha(x)+\beta}}{c(x)} \biggl(\int_{\R
}p(x,\mathrm{d}y)V(y)-V(x)
\biggr)>0.
\]
We have
%
\begin{eqnarray}\label{eq40}
&&\int_{\R}p(x,\mathrm{d}y)V(y)-V(x)\nonumber\\
&&\quad=\int_{\R
}V(y+x)f_x(y)\,\mathrm{d}y-V(x)
\\
&&\quad=\int_{-x}^{\infty} \bigl(1-(1+y+x)^{-\beta}
\bigr)f_x(y)\,\mathrm{d}y+\int_{-\infty}^{-x}
\bigl(1-(1-y-x)^{-\beta} \bigr)f_x(y)\,\mathrm{d}y
\nonumber
\\
&&\qquad{} - \bigl(1-\bigl(1+|x|\bigr)^{-\beta} \bigr)\int_{-x}^{\infty
}f_x(y)\,\mathrm{d}y-
\bigl(1-\bigl(1+|x|\bigr)^{-\beta} \bigr)\int_{-\infty
}^{-x}f_x(y)\,\mathrm{d}y\nonumber
\\
&&\quad=\bigl(1+|x|\bigr)^{-\beta} \biggl[\int_{-x}^{\infty}
\biggl(1- \biggl(\frac{1+|x|}{1+x+y} \biggr)^{\beta} \biggr)f_x(y)\,\mathrm{d}y
\nonumber
\\
&&\hspace*{70pt}{} +\int_{-\infty}^{-x} \biggl(1- \biggl(
\frac
{1+|x|}{1-x-y} \biggr)^{\beta} \biggr)f_x(y)\,\mathrm{d}y \biggr].\nonumber
\end{eqnarray}

\textit{Step} 2. In the second step, by use of condition (C3), we
find an operable lower bound for (\ref{eq40}). First, let us take
a look at the case when $x>0$. Let $0<\varepsilon<1$ be arbitrary.
Then, by (C3), there exists $y_{\varepsilon}\geq a_0\vee1$ (the
constant $a_0>0$ is defined in (\ref{eq13})), such that for all
$|y|\geq y_{\varepsilon}$
\[
\biggl\llvert f_x(y)\frac{|y|^{\alpha(x)+1}}{c(x)}-1\biggr\rrvert <\varepsilon
\]
for all $x\in[-k,k]^{c}.$
Let $x\geq k\vee y_{\varepsilon}.$ Then we have
\begin{eqnarray*}
\int_{-x}^{\infty} \biggl(1- \biggl(1+
\frac
{y}{1+x} \biggr)^{-\beta} \biggr)f_x(y)\,\mathrm{d}y&>& c(x) (1+
\varepsilon)\int_{y_{\varepsilon}}^{x} \biggl(1- \biggl(1-
\frac{y}{1+x} \biggr)^{-\beta} \biggr)\frac{\mathrm{d}y}{y^{\alpha
(x)+1}}
\\
&&{} +\int_{-y_{\varepsilon}}^{y_{\varepsilon}} \biggl(1- \biggl(1+
\frac{y}{1+x} \biggr)^{-\beta} \biggr)f_x(y)\,\mathrm{d}y
\\
&&{} +c(x) (1-\varepsilon)\int_{y_{\varepsilon}}^{\infty} \biggl(1-
\biggl(1+\frac{y}{1+x} \biggr)^{-\beta} \biggr)\frac
{\mathrm{d}y}{y^{\alpha(x)+1}}
\end{eqnarray*}
and
\begin{eqnarray*}
&&\int_{-\infty}^{-x} \biggl(1- \biggl(
\frac
{1+x}{1-x-y} \biggr)^{\beta} \biggr)f_x(y)\,\mathrm{d}y
\\
&&\quad> \frac{c(x)(1-\varepsilon)}{\alpha(x)x^{\alpha
(x)}}-c(x) (1+\varepsilon )\int_{x}^{\infty}
\biggl(\frac{1+x}{1-x+y} \biggr)^{\beta
}\frac{\mathrm{d}y}{y^{\alpha(x)+1}}.
\end{eqnarray*}
Note that this was the crucial step where we needed condition (C3).
For given $0<\varepsilon<1$ and $x\geq k\vee y_{\varepsilon}$, let us
put
\begin{eqnarray*}
U^{\varepsilon}_1(x)&:=&c(x) (1+\varepsilon)\int
_{y_{\varepsilon}}^{x} \biggl(1- \biggl(1-\frac{y}{1+x}
\biggr)^{-\beta
} \biggr)\frac{\mathrm{d}y}{y^{\alpha(x)+1}},
\\
U^{\varepsilon}_2(x)&:=&\int_{-y_{\varepsilon}}^{y_{\varepsilon
}}
\biggl(1- \biggl(1+\frac{y}{1+x} \biggr)^{-\beta}
\biggr)f_x(y)\,\mathrm{d}y,
\\
U^{\varepsilon}_3(x)&:=&c(x) (1-\varepsilon)\int
_{y_{\varepsilon}}^{\infty} \biggl(1- \biggl(1+\frac{y}{1+x}
\biggr)^{-\beta
} \biggr)\frac{\mathrm{d}y}{y^{\alpha(x)+1}},
\\
U^{\varepsilon}_4(x)&:=&\frac{c(x)(1-\varepsilon)}{\alpha
(x)x^{\alpha(x)}} \quad\mbox{and}
\\
U^{\varepsilon}_5(x)&:=&c(x) (1+\varepsilon)\int
_{x}^{\infty
} \biggl(\frac{1+x}{1-x+y}
\biggr)^{\beta}\frac{\mathrm{d}y}{y^{\alpha(x)+1}}.
\end{eqnarray*}
Hence, we have
\[
\int_{\R}p(x,\mathrm{d}y)V(y)-V(x)>U^{\varepsilon}_1(x)+U^{\varepsilon
}_2(x)+U^{\varepsilon}_3(x)+U^{\varepsilon}_4(x)-U^{\varepsilon}_5(x)
.
\]

\textit{Step} 3. In the third step, we prove
%
\begin{eqnarray}\label{eq417}
&&\hspace*{-10pt}\operatorname{\lim\inf}\limits_{x\longrightarrow\infty}\frac
{\alpha(x)x^{\alpha(x)+\beta}}{c(x)} \biggl(\int
_{\R
}p(x,\mathrm{d}y)V(y)-V(x) \biggr)
\nonumber
\\[-8pt]\\[-8pt]
&&\hspace*{-10pt}\quad>\operatorname{\lim\inf}\limits_{\varepsilon\longrightarrow0}\operatorname{\lim\inf}\limits_{y_{\varepsilon
}\longrightarrow\infty}\operatorname{\lim\inf}\limits_{x\longrightarrow\infty}
\frac
{\alpha
(x)x^{\alpha(x)}}{c(x)} \bigl(U^{\varepsilon}_1(x)+U^{\varepsilon
}_3(x)+U^{\varepsilon}_4(x)-U^{\varepsilon}_5(x)
\bigr)-T(\alpha ,\beta )\geq0.\quad\nonumber
\end{eqnarray}
Recall that
\[
T(\alpha,\beta)={}_2F_1(-\alpha,\beta,1-y;1)+\beta B(1;
\alpha+\beta,1-\alpha)-\alpha B(1;\alpha+\beta,1-\beta)
\]
and
\[
\operatorname{\lim\inf}\limits_{\varepsilon\longrightarrow0}\operatorname{\lim\inf}\limits_{y_{\varepsilon
}\longrightarrow\infty}\operatorname{\lim\inf}\limits_{x\longrightarrow\infty}
\frac
{\alpha
(x)x^{\alpha(x)}}{c(x)}U^{\varepsilon}_2(x)>-T(\alpha,\beta)
\]
(assumption (\ref{eq14})). By straightforward calculations, using
(\ref{eq42}), (\ref{eq47}) and (\ref{eq43}), we have
\begin{eqnarray*}
\frac{\alpha(x)x^{\alpha(x)}}{c(x)}U^{\varepsilon
}_1(x)&=&\frac{(1+\varepsilon)x^{\alpha(x)}}{y_{\varepsilon}^{\alpha(x)}} -
\frac{(1+\varepsilon)x^{\alpha(x)}{}_2F_1 (-\alpha(x),\beta,
1 - \alpha(x); {y_{\varepsilon}}/({1 +
x}) )}{y_{\varepsilon}^{\alpha(x)}}
\\
&&{}-(1+\varepsilon)+(1+\varepsilon){}_2F_1 \biggl(-
\alpha(x),\beta, 1 - \alpha(x); \frac{x}{1+x} \biggr),
\\
\frac{\alpha(x)x^{\alpha(x)}}{c(x)}U^{\varepsilon
}_4(x)&=&(1-\varepsilon )\quad
\mbox{and}
\\
\frac{\alpha(x)x^{\alpha(x)}}{c(x)}U^{\varepsilon}_5(x)&=& \frac{(1+\varepsilon)\alpha(x)x^{\alpha(x)}(1+x)^{\beta}B
(({x-1})/{x};
\alpha(x) + \beta,
1-\beta )}{(x-1)^{\alpha(x)+\beta}}.
\end{eqnarray*}
It is easy to
check that
\[
\frac{\partial}{\partial y} \biggl(-\frac{{}_2F_1 (-\alpha
(x),\beta
,1-\alpha(x);-{y}/({1+x}) )}{\alpha(x)y^{\alpha
(x)}(1+x)^{\beta
}} \biggr)=\frac{1}{(1+x+y)^{\beta}y^{\alpha(x)+1}},
\]
and from (\ref{eq48}) and
%
\begin{equation}
\label{eq410} \Gamma(z+1)=z\Gamma(z),\qquad z\in\CC\setminus\ZZ_-,
\end{equation}
we have
\begin{eqnarray*}
\frac{{}_2F_1 (-\alpha(x),\beta,1-\alpha(x);-
{y}/({1+x}) )}{\alpha(x)y^{\alpha(x)}(1+x)^{\beta}}&=&\frac
{{}_2F_1
(\beta,\alpha(x)+\beta,1+\alpha(x)+\beta;-({1+x})/{y}
)}{(\alpha
(x)+\beta)y^{\alpha(x)+\beta}}
\\
&&{}+\frac{\Gamma(1-\alpha(x))\Gamma
(\alpha
(x)+\beta)}{\alpha(x)(1+x)^{2\alpha(x)+\beta}\Gamma(\beta)}.
\end{eqnarray*}
Therefore,
\[
\int\frac{\mathrm{d}y}{(1+x+y)^{\beta}y^{\alpha(x)+1}}=-\frac{{}_2F_1
(\beta
,\alpha(x)+\beta,1+\alpha(x)+\beta;-({1+x})/{y} )}{(\alpha
(x)+\beta)y^{\alpha(x)+\beta}},
\]
that is,
\begin{eqnarray*}
&&\frac{\alpha(x)x^{\alpha(x)}}{c(x)}U^{\varepsilon
}_3(x)
\\
&&\quad=\frac{(1-\varepsilon)x^{\alpha(x)}}{y_{\varepsilon}^{\alpha
(x)}}-\frac
{(1-\varepsilon)\alpha(x)x^{\alpha(x)}(1+x)^{\beta}{}_2F_1 (\beta,\alpha(x)+\beta,1+\alpha(x)+\beta,-
({1+x})/{y_{\varepsilon}} )}{y_{\varepsilon}^{\alpha(x)+\beta
}(\alpha
(x)+\beta)}.
\end{eqnarray*}
Furthermore, from (\ref{eq45}), (\ref{eq48}) and (\ref{eq410}),
we have
\begin{eqnarray*}
&&\frac{\alpha(x)x^{\alpha(x)}}{c(x)}U^{\varepsilon
}_3(x)
\\
&&\quad=\frac{(1-\varepsilon)x^{\alpha(x)}}{y_{\varepsilon}^{\alpha(x)}} - \frac{(1-\varepsilon)x^{\alpha(x)}}{y_{\varepsilon}^{\alpha(x)}}{}_2F_1
\biggl(\beta,-\alpha(x),1-\alpha(x);-\frac{y_{\varepsilon
}}{x+1} \biggr)
\\
&&\qquad {}-\frac{(1-\varepsilon)\Gamma(\alpha(x)+\beta)\Gamma(-\alpha
(x))\alpha(x)x^{\alpha(x)}}{\Gamma(\beta)(1+x)^{\alpha(x)}}{}_2F_1 \biggl(\alpha(x)+
\beta,0,1+\alpha(x);-\frac{y_{\varepsilon}}{x+1} \biggr)
\\
&&\quad=\frac{(1-\varepsilon)x^{\alpha(x)}}{y_{\varepsilon}^{\alpha(x)}} - \frac{(1-\varepsilon)x^{\alpha(x)}}{y_{\varepsilon}^{\alpha(x)}}{}_2F_1
\biggl(\beta,-\alpha(x),1-\alpha(x);-\frac{y_{\varepsilon
}}{x+1} \biggr)
\\
&&\qquad{} -\frac{(1-\varepsilon)\Gamma(\alpha(x)+\beta)\Gamma(-\alpha
(x))\alpha
(x)x^{\alpha(x)}}{\Gamma(\beta)(1+x)^{\alpha(x)}}.
\end{eqnarray*}
Let us put
\begin{eqnarray*}
V_1^{\varepsilon}(x)&:=&\frac{(1+\varepsilon)x^{\alpha
(x)}}{y_{\varepsilon}^{\alpha(x)}}- \frac{(1+\varepsilon)x^{\alpha(x)}{}_2F_1 (-\alpha(x),\beta, 1
- \alpha(x); {y_{\varepsilon}}/({1 +
x}) )}{y_{\varepsilon}^{\alpha(x)}}
\\
V_2^{\varepsilon}(x)&:=&\frac{(1-\varepsilon)x^{\alpha
(x)}}{y_{\varepsilon}^{\alpha(x)}} - \frac{(1-\varepsilon)x^{\alpha(x)}}{y_{\varepsilon}^{\alpha(x)}}{}_2F_1
\biggl(\beta,-\alpha(x),1-\alpha(x);-\frac{y_{\varepsilon
}}{x+1} \biggr)
\end{eqnarray*}
and
\begin{eqnarray*}
V_3^{\varepsilon}(x)&:=&(1+\varepsilon){}_2F_1
\biggl(-\alpha(x),\beta , 1 - \alpha(x); \frac{x}{1+x} \biggr)
\\
&&{}-\frac{(1-\varepsilon)\Gamma(\alpha
(x)+\beta
)\Gamma(-\alpha(x))\alpha(x)x^{\alpha(x)}}{\Gamma(\beta
)(1+x)^{\alpha
(x)}}
\\
&&{}-\frac{(1+\varepsilon)\alpha(x)x^{\alpha(x)}(1+x)^{\beta
}B
(({x-1})/{x};
\alpha(x) + \beta,
1-\beta )}{(x-1)^{\alpha(x)+\beta}}.
\end{eqnarray*}
Hence,
(\ref{eq417}) is reduced to
%
\begin{eqnarray}
\label{eq411}
&&\operatorname{\lim\inf}\limits_{x\longrightarrow\infty}\frac
{\alpha(x)x^{\alpha(x)+\beta}}{c(x)} \biggl(\int
_{\R
}p(x,\mathrm{d}y)V(y)-V(x) \biggr)
\nonumber
\\
&&\quad>\operatorname{\lim\inf}\limits_{\varepsilon
\longrightarrow
0}\operatorname{\lim\inf}\limits_{y_{\varepsilon}\longrightarrow\infty}\operatorname{\lim\inf}\limits_{x\longrightarrow\infty}V^{\varepsilon}_1(x)+
\operatorname{\lim\inf}\limits_{\varepsilon
\longrightarrow0}\operatorname{\lim\inf}\limits_{y_{\varepsilon}\longrightarrow\infty
}\operatorname{\lim\inf}\limits_{x\longrightarrow\infty}V_2^{\varepsilon}(x)
\\
&&\qquad {}+\operatorname{\lim\inf}\limits_{\varepsilon\longrightarrow0}\operatorname{\lim\inf}\limits_{x\longrightarrow
\infty
}V_3^{\varepsilon}(x)-T(
\alpha,\beta).\nonumber
\end{eqnarray}

By (\ref{eq41}) and (\ref{eq42}), we have
\[
0\leq{}_2F_1 \biggl(-\alpha(x),\beta, 1 - \alpha(x),
\frac{y_{\varepsilon}}{1
+ x} \biggr)\leq1,
\]
therefore
%
\begin{equation}
\label{eq412}\operatorname{\lim\inf}\limits_{\varepsilon\longrightarrow
0}\operatorname{\lim\inf}\limits_{y_{\varepsilon}\longrightarrow\infty}\operatorname{\lim\inf}\limits_{x\longrightarrow
\infty
}V^{\varepsilon}_1(x)
\geq0.
\end{equation}

Since $1-\alpha(x)-(-\alpha(x))-\beta=1-\beta>0$, from
(\ref{eq41}) and the dominated convergence theorem, we have
%
\begin{equation}
\label{eq413}\operatorname{\lim\inf}\limits_{\varepsilon\longrightarrow
0}\operatorname{\lim\inf}\limits_{y_{\varepsilon}\longrightarrow\infty}\operatorname{\lim\inf}\limits_{x\longrightarrow
\infty
}V^{\varepsilon}_2(x)=0.
\end{equation}

At the end, let us calculate
\[
\operatorname{\lim\inf}\limits_{\varepsilon\longrightarrow0}\operatorname{\lim\inf}\limits_{x\longrightarrow
\infty
}V_3^{\varepsilon}(x).
\]
From (\ref{eq41}), we have
\[
{}_2F_1 \biggl(-\alpha(x), \beta, 1 - \alpha(x);
\frac{x}{1+x} \biggr)\geq{}_2F_1 \bigl(-\alpha(x),
\beta, 1 - \alpha(x); 1 \bigr),
\]
and from (\ref{eq43}) we have
\[
\alpha(x)B \biggl(\frac{x-1}{x}; \alpha(x) + \beta, 1 - \beta \biggr)\leq
\alpha(x)B \bigl(1; \alpha(x) + \beta, 1 - \beta \bigr).
\]
Hence, we have
\begin{eqnarray*}
&&\operatorname{\lim\inf}\limits_{\varepsilon\longrightarrow0}\operatorname{\lim\inf}\limits_{x\longrightarrow\infty}V_3^{\varepsilon}(x)
\\
&&\quad\geq\operatorname{\lim\inf}\limits_{\varepsilon
\longrightarrow0}\operatorname{\lim\inf}\limits_{x\longrightarrow\infty} \biggl[(1+
\varepsilon){}_2F_1\bigl(-\alpha(x),\beta,1-\alpha(x);1
\bigr)
\nonumber
\\
&&\hspace*{77pt}{}-\frac
{(1-\varepsilon)\alpha(x)\Gamma(\alpha(x)+\beta)\Gamma(-\alpha
(x))x^{\alpha(x)}}{\Gamma(\beta)(1+x)^{\alpha(x)}}
\nonumber
\\
&&\hspace*{77pt}{}-\frac{(1+\varepsilon)\alpha(x)B(1;\alpha(x)+\beta,1-\beta
)x^{\alpha
(x)}(1+x)^{\beta}}{(1-x)^{\alpha(x)+\beta}} \biggr],
\end{eqnarray*}
that is, since all terms are bounded,
\begin{eqnarray*}
&&\operatorname{\lim\inf}\limits_{\varepsilon\longrightarrow0}\operatorname{\lim\inf}\limits_{x\longrightarrow\infty}V_3^{\varepsilon}(x)
\\
&&\quad\geq\operatorname{\lim\inf}\limits_{x\longrightarrow\infty} \bigl[{}_2F_1\bigl(-\alpha
(x),\beta ,1-\alpha(x);1\bigr)+\beta B\bigl(1;\alpha(x)+\beta,1-\alpha(x)\bigr)
\\
&&\hspace*{48pt}{}-\alpha(x)B\bigl(1;\alpha(x)+\beta ,1-\beta\bigr) \bigr].
\end{eqnarray*}
One can
prove that the function
\[
y\longmapsto T(y,\beta):=\label{eq414}{}_2F_1(-y,
\beta,1-y;1)+\beta B(1;y+\beta,1-y)-y B(1;y+\beta,1-\beta)
\]
is strictly decreasing on
$[0,1-\beta)$, and it easy to see that $T(1-\beta,\beta)=0.$ Hence,
since $0\leq\alpha<1-\beta$, we have
%
\begin{equation}
\label{eq415}\operatorname{\lim\inf}\limits_{\varepsilon\longrightarrow
0}\operatorname{\lim\inf}\limits_{x\longrightarrow\infty}V_3^{\varepsilon}(x)
\geq T(\alpha,\beta).
\end{equation}
By combining (\ref{eq411}), (\ref{eq412}),
(\ref{eq413}) and (\ref{eq415}), we have
\[
\operatorname{\lim\inf}\limits_{x\longrightarrow\infty}\frac{\alpha(x)x^{\alpha(x)+\beta
}}{c(x)} \biggl(\int_{\R}p(x,\mathrm{d}y)V(y)-V(x)
\biggr)>0.
\]
The case when $x<0$ is treated in the same way. Therefore, by
\cite{meyn}, Theorem 8.4.3, the chain $\{X_n\}_{n\geq0}$ is
transient.
\end{pf*}

\section{Some remarks and generalizations of the model}\label{sec5}
We start this section with the proof of equivalence of
conditions (\ref{eq13}) and (\ref{eq15}), and the proof of
relaxation of condition (\ref{eq14}) to condition
(\ref{eq17}).
\begin{longlist}[(ii)]
\item[(i)]
Recall that condition (\ref{eq15}) is given by
\[
\operatorname{\lim\sup}\limits_{\delta\longrightarrow0}\operatorname{\lim\sup}\limits_{|x|\longrightarrow\infty}\frac{(1+|x|)^{\alpha(x)}}{c(x)}\int
_{-\delta
(1+|x|)}^{\delta(1+|x|)}\ln \biggl(1+\sgn (x)
\frac{y}{1+|x|} \biggr)f_{x}(y)\,\mathrm{d}y<R(\alpha).
\]
Using $\ln(1+t)\leq t$, condition (\ref{eq15}) follows from the
condition
%
\begin{equation}
\label{eq51}\operatorname{\lim\sup}\limits_{\delta\longrightarrow0}\operatorname{\lim\sup}\limits_{|x|\longrightarrow\infty}\sgn(x)
\frac{(1+|x|)^{\alpha(x)-1}}{c(x)}\int_{-\delta
(1+|x|)}^{\delta(1+|x|)}y
f_{x}(y)\,\mathrm{d}y<R(\alpha).
\end{equation}
In fact, under condition (C3), conditions (\ref{eq15}) and
(\ref{eq51}) are equivalent, but the proof of this statement is
rather elementary and technical and we omit it here. Furthermore, by
(C3) and since $\alpha(x)\in(1,2)$, condition (\ref{eq51}) is
equivalent with
\[
\operatorname{\lim\sup}\limits_{|x|\longrightarrow\infty}\sgn(x)\frac{(1+|x|)^{\alpha(x)-1}}{c(x)}\int
_{\R}y f_{x}(y)\,\mathrm{d}y<R(\alpha),
\]
that is, with condition (\ref{eq13}). Indeed, let $\delta>0$ and
$0<\varepsilon<1$ be arbitrary. Then, by (C3), there exists
$y_\varepsilon>0$ such that for all $|y|\geq y_{\varepsilon}$
\[
\biggl\llvert f_x(y)\frac{|y|^{\alpha(x)+1}}{c(x)}-1\biggr\rrvert <\varepsilon
\]
for all $x\in[-k,k]^{c}.$ By taking
$|x|\geq\frac{y_{\varepsilon}}{\delta}-1,$ we have (recall that
$\alpha(x)\in(1,2)$)
\begin{eqnarray*}
\int_{\R}y f_x(y)\,\mathrm{d}y &>&-(1+\varepsilon)\int
_{-\infty}^{-\delta(1+|x|)}\frac
{c(x)}{|y|^{\alpha
(x)}}\, \mathrm{d}y+\int
_{-\delta(1+|x|)}^{\delta
(1+|x|)}y f_x(y)\,\mathrm{d}y
\\
&&{}+(1-\varepsilon)\int_{\delta(1+|x|)}^{\infty
}
\frac
{c(x)}{|y|^{\alpha(x)}}\,\mathrm{d}y
\\
&=&\int_{-\delta(1+|x|)}^{\delta
(1+|x|)}y f_x(y)\,\mathrm{d}y-
\frac{2\varepsilon
c(x)}{(\alpha(x)-1)\delta^{\alpha(x)-1}(1+|x|)^{\alpha(x)-1}}.
\end{eqnarray*}
In the same way, we get
\[
\int_{\R}y f_x(y)\,\mathrm{d}y<\int_{-\delta(1+|x|)}^{\delta
(1+|x|)}y
f_x(y)\,\mathrm{d}y+\frac{2\varepsilon
c(x)}{(\alpha(x)-1)\delta^{\alpha(x)-1}(1+|x|)^{\alpha(x)-1}}.
\]
By taking
$\operatorname{\lim\sup}_{|x|\longrightarrow\infty}$,
$\operatorname{\lim\sup}_{\varepsilon\longrightarrow0}$ and
$\operatorname{\lim\sup}_{\delta\longrightarrow0}$ we get the desired result.
\item[(ii)] From the concavity of the function $x\longmapsto x^{\beta
}$, for
$\beta\in(0,1-\alpha)$, we have
\begin{eqnarray*}
&&\operatorname{\lim\inf}\limits_{|x|\longrightarrow\infty}\frac
{\alpha(x)|x|^{\alpha(x)}}{c(x)}\int_{-a}^{a}
\biggl(1- \biggl(1+\sgn(x)\frac{y}{1+| x|} \biggr)^{-\beta}
\biggr)f_x(y)\,\mathrm{d}y
\\
&&\quad\geq\operatorname{\lim\inf}\limits_{|x|\longrightarrow\infty}\frac{\alpha(x)|x|^{\alpha
(x)}}{c(x)}\frac{(1+|x|-a)^{\beta}-(1+|x|)^{\beta}}{(1+|x|-a)^{\beta
}}\int
_{-a}^{a}f_x(y)\,\mathrm{d}y
\\
&&\quad\geq\operatorname{\lim\inf}\limits_{|x|\longrightarrow\infty} \biggl(-\frac{a\beta\alpha
(x)|x|^{\alpha(x)}}{c(x)(1+|x|-a)} \biggr)=-a\beta
\operatorname{\lim\sup}\limits_{|x|\longrightarrow\infty}\frac{\alpha(x)}{c(x)}|x|^{\alpha(x)-1}.
\end{eqnarray*}
\end{longlist}

In the sequel, we give several generalizations of the stable-like
chain $\{X_n\}_{n\geq0}$. Recall that a function
$f\dvtx \R\longrightarrow\R$ is called \emph{lower semicontinuous} if
$\operatorname{\lim\inf}_{y\longrightarrow x}f(y)\geq f(x)$ for all $x\in\R$.
%
\begin{definition} Let $\{Y_n\}_{n\geq0}$ be a Markov chain on
$(\R,\mathcal{B}(\R)).$
\begin{enumerate}[(iii)]
\item[(i)]The chain $\{Y_n\}_{n\geq0}$ is called a \emph{T-chain} if
for some probability measure $a=\{a(n)\}_{n\geq0}$ on $\ZZ_+$ there
exists a kernel $T(x,B)$ on $(\R,\mathcal{B}(\R))$ with
$T(x,\R) > 0$ for all $x\in\R$, such that the function\vadjust{\goodbreak} $x\longmapsto
T(x,B)$ is lower semicontinuous for all $B\in\mathcal{B}(\R)$, and
\[
\sum_{n=0}^{\infty}a(n)p^{n}(x,B)
\geq T(x,B)
\]
holds for all
$x\in\R$ and all $B\in\mathcal{B}(\R).$
\item[(ii)]The chain $\{Y_n\}_{n\geq0}$ is \emph{Harris
recurrent}, or \emph{H-recurrent},
if it is
$\psi$-irreducible and if $\mathbb{P}(\tau_B<\infty|Y_0=x)=1$ holds
for all $x\in\R$ and all $B\in\mathcal{B}^{+}(\R)$.
\item[(iii)]A state $x\in\R$ is called a \emph{topologically
recurrent state}
if $\sum_{n=0}^{\infty}p^{n}(x,O_x)=\infty$ holds for
all open neighborhoods $O_x$ around $x$. Otherwise we call state $x$ a
\emph{topologically transient state.}
\end{enumerate}
\end{definition}
From Proposition \ref{p26} and \cite{meyn}, Theorem 6.2.5, we
have the following.
%
\begin{proposition} The chain $\{X_n\}_{n\geq0}$ is a T-chain.
\end{proposition}

It is well known that the
recurrence and H-recurrence properties of a Markov chain on the
general state space are not equivalent (see \cite{meyn}, Section 9.1.2).
Now, let us prove that these properties are equivalent for the
stable-like chain $\{X_n\}_{n\geq0}$.
%
\begin{proposition} The chain $\{X_n\}_{n\geq0}$ is recurrent if and
only if it is H-recurrent.
\end{proposition}
\begin{pf}
We have to prove that recurrence property implies H-recurrence
property, since the opposite claim is trivial. Since the Markov
chain $\{X_n\}_{n\geq0}$ is a T-chain, by \cite{meyn}, Theorem
9.3.6, it is enough to prove that every state is a
topologically recurrent state. That follows from \cite{meyn}, Lemma~6.1.4 and Theorem 9.3.3.
\end{pf}

If we change the chain $\{X_n\}_{n\geq0}$ on a set of Lebesgue
measure zero, it can happen that its recurrence and H-recurrence
properties are not equivalent anymore. Let $A\in\mathcal{B}(\R)$
be such that $\lambda(A)=0$. Note that $A$ can be unbounded. Let
$\{\bar{X}_n\}_{n\geq0}$ be a Markov chain on
$(\R,\mathcal{B}(\R))$ given by the transition kernel
\[
\bar{p}(x,\mathrm{d}y)=\bar{f}_x(y-x)\,\mathrm{d}y,
\]
where
$\{\bar{f}_x\dvt  x\in\R\}$ is the family of density functions on $\R$
such that $\bar{f}_x=f_x$, for every $x\in\R\setminus A$. It is to
easy see that the chain $(\bar{X}_n)$ is $\lambda$-irreducible and
aperiodic. Therefore, a Borel
set is a small set for $\{\bar{X}_n\}_{n\geq0}$ if and only if it is
a petite set for $\{\bar{X}_n\}_{n\geq0}$. But we cannot conclude
that every bounded Borel set is a petite set. The most we can get is
that every bounded set $B\in\mathcal{B}(\R\setminus A)$ is a petite
set. As a consequence of this fact, we do not know if the chain
$\{\bar{X}_n\}_{n\geq0}$ is a T-chain, so we cannot deduce
equivalence between recurrence and H-recurrence property of the
chain $\{\bar{X}_n\}_{n\geq0}$. But, since the chains
$\{X_n\}_{n\geq0}$ and $\{\bar{X}_n\}_{n\geq0}$ are
$\lambda$-irreducible and since they differ on the set with zero
Lebesgue measure, it is easy to see that the recurrence property of
the chain $\{X_n\}_{n\geq0}$ is equivalent with the recurrence
property of the chain $\{\bar{X}_n\}_{n\geq0}$, and the H-recurrence
property of the chain $\{X_n\}_{n\geq0}$ is equivalent with the
H-recurrence property of the chain $\{\bar{X}_n\}_{n\geq0}$. Hence,
the chain $\{\bar{X}_n\}_{n\geq0}$ is recurrent if and only if it is
H-recurrent.

In Proposition \ref{p26}, it is proved that every bounded Borel set
is a petite set (singleton) for the stable-like chain
$\{X_n\}_{n\geq0}$. Therefore, it is natural to expect that a change
of the chain $\{X_n\}_{n\geq0}$ on an arbitrary bounded Borel set
will not affect its recurrence and transience property. Let
$B\in\mathcal{B}(\R)$ be bounded and let $\{\tilde{X}_n\}_{n\geq0}$
be a stable-like Markov chain on $(\R,\mathcal{B}(\R))$ given by
the transition kernel
\[
\tilde{p}(x,\mathrm{d}y)=\tilde{f}_x(y-x)\,\mathrm{d}y,
\]
where
$\{\tilde{f}_x\dvt  x\in\R\}$ is a family of density functions on $\R$
such that $\tilde{f}_x=f_x$ for all $x\in\R\setminus B$ and such
that it satisfies conditions (C1)--(C5). Therefore, the chain
$\{\tilde{X}_n\}_{n\geq0}$ is either H-recurrent or transient.

\begin{proposition} The chain $\{X_n\}_{n\geq0}$ is H-recurrent if and
only if the chain $\{\tilde{X}_n\}_{n\geq0}$ is
H-recurrent. Hence, the chain $\{X_n\}_{n\geq0}$ is recurrent if and
only if the chain $\{\tilde{X}_n\}_{n\geq0}$ is recurrent.
\end{proposition}
\begin{pf} If $\lambda(B)=0$, the claim follows from the above discussion.
Let us suppose that $\lambda(B)>0$. By Proposition \ref{p26}, the
set $B$ is a petite set for both chains $\{X_n\}_{n\geq0}$ and
$\{\tilde{X}_n\}_{n\geq0}$. Let us suppose that the chain
$\{X_n\}_{n\geq0}$ is H-recurrent. Then, by \cite{meyn}, Theorem
9.1.4, we have $\mathbb{P}(\tau_B<\infty|X_0=x)=1$
for all $x\in\R$. Since
\[
\mathbb{P}(\tau_B<\infty|X_0=y)=\mathbb{P}(\tilde{
\tau}_B<\infty |\tilde{X}_0=y)
\]
for all $y\notin B$,
we have
\begin{eqnarray*}
\mathbb{P}(\tilde{\tau}_B<\infty|\tilde {X}_0=x)&=&
\tilde {p}(x,B)+\int_{B^{c}}\hat{p}(x,\mathrm{d}y)\mathbb{P}(\tilde{\tau
}_B<\infty|\tilde{X}_0=y)
\\
&=&\tilde{p}(x,B)+\tilde{p}\bigl(x,B^{c}\bigr)=1
\end{eqnarray*}
for all $x\in\R$. Therefore, by \cite{meyn}, Proposition 9.1.7,
the chain $\{\tilde{X}_n\}_{n\geq0}$ is H-recurrent. The proof of
the opposite direction is completely the same.
\end{pf}

From the above discussions, we can weaken assumptions on function
$\alpha(x)$ and conditions (\ref{eq13}) and (\ref{eq14}) in
Theorems \ref{tm13} and \ref{tm14}.
In Theorem \ref{tm13}, we assumed that
$\alpha\dvtx \R\longrightarrow(1,2)$ and
\[
\operatorname{\lim\inf}\limits_{|x|\longrightarrow\infty}\alpha(x)>1,
\]
but it is enough to request that $\alpha:\R\setminus(A\cup
B)\longrightarrow(1,2)$
and
\[
\operatorname{\lim\inf}\limits_{x\in\R\setminus A,
|x|\longrightarrow\infty}\alpha(x)>1
\]
for some set
$A\in\mathcal{B}(\R)$ with zero Lebesgue measure and some bounded
set $B\in\mathcal{B}(\R)$. In condition (\ref{eq13}) instead of
using $\operatorname{\lim\sup}_{|x|\longrightarrow\infty}$, we use $\operatorname{\lim\sup}_{x\in
\R\setminus A, |x|\longrightarrow\infty}$. An analog modification
can be done in Theorem \ref{tm14}.

The transition densities of the stable-like chain
$\{X_n\}_{n\geq0}$, from the current state $x$, have the power-law
decay with exponent $\alpha(x)+1$. Let us take a look at the Markov
chain with transition densities with the power-law decay with
exponent $\alpha_-(x)+1$ on the left of the current state $x$
and with the
power-law decay with exponent $\alpha_+(x)+1$ on the right of the
current state $x$. Let $\alpha_+, \alpha_-\dvtx \R\longrightarrow(0,2)$
and $c_+, c_-\dvtx \R\longrightarrow(0,\infty)$ be arbitrary functions
and let $(X'_n)$ be a Markov chain on $(\R,\mathcal{B}(\R))$ given
by the transition kernel $p'(x,\mathrm{d}y)=f'_x(y-x)\,\mathrm{d}y$, where
$\{f'_x\dvtx x\in\R\}$ is a family of density functions on $\R$ which
satisfies:
\begin{enumerate}[(C1$^\prime$)]
\item[(C1$^\prime$)] $x\longmapsto f'_x(y)$ is measurable, for every
$y\in\R$;
\item[(C2$^\prime$)] $f'_x(y)\sim c_+(x)y^{-\alpha_+(x)-1},$
when
$y\longrightarrow\infty$, and $f'_x(y)\sim c_-(x)(-y)^{-\alpha_-(x)-1},$
when
$y\longrightarrow-\infty$;
\item[(C3$^\prime$)]there exists $k'>0$ such that
\[
\lim_{y\longrightarrow\infty}\sup_{x\in[-k',k']^{c}}\biggl\llvert f'_x(y)
\frac
{y^{\alpha_+(x)+1}}{c_+(x)}-1\biggr\rrvert =0
\]
and
\[
\lim_{y\longrightarrow-\infty}\sup_{x\in[-k',k']^{c}}\biggl\llvert f'_x(y)
\frac
{(-y)^{\alpha_-(x)+1}}{c_-(x)}-1\biggr\rrvert =0;
\]
\item[(C4$^\prime$)] $\inf_{x\in C}(c_+(x)\wedge c_-(x))>0$
for every compact set
$C\subseteq[-k',k']^{c}$;
\item[(C5$^\prime$)] there exists $l'>0$ such that for every compact
set $C\subseteq[-l',l']^{c}$ with $\lambda(C)>0$,
we have
\[
\inf_{x\in[-k',k']}\int_{C-x}f'_x(y)\,\mathrm{d}y>0.
\]
\end{enumerate}
It is clear that the chain $\{X'_n\}_{n\geq0}$ has the same
properties, discussed in Section \ref{sec2}, as the chain $\{X_n\}_{n\geq0}$.
It is $\lambda$-irreducible and aperiodic and every bounded Borel
set is a petite set. By assuming certain additional conditions,
Theorems \ref{tm13} and \ref{tm14} can be generalized in terms of
the chain $\{X'_n\}_{n\geq0}$. The chain
$\{X'_n\}_{n\geq0}$ will be recurrent if $\alpha_+, \alpha_-\dvtx \R
\longrightarrow(1,2)$ are such that
\[
\lim_{|x|\longrightarrow\infty}\frac{\alpha_+(x)}{\alpha_-(x)}=1 \quad\mbox{and}\quad\alpha:=
\operatorname{\lim\inf}\limits_{|x|\longrightarrow\infty
}\alpha_+(x) \Bigl(=\operatorname{\lim\inf}\limits_{|x|\longrightarrow\infty}\alpha_-(x)
\Bigr)>1,
\]
and
$c_+, c_-\dvtx \R\longrightarrow(0,\infty)$ are such that
\[
\lim_{|x|\longrightarrow\infty}\frac{c_-(x)}{c_+(x)}|x|^{\alpha
_+(x)-\alpha_-(x)}=1
\]
and such that condition (\ref{eq13}) is satisfied with the
constant $R(\alpha)$. In this case, for the test function $V(x)$, we
take $V(x)=\ln(1+|x|)$ again.
Similarly, the chain
$(X'_n)$ will be transient if $\alpha_+, \alpha_-\dvtx \R\longrightarrow
(0,1)$ are such that
\[
\alpha_+:=\operatorname{\lim\sup}\limits_{|x|\longrightarrow\infty}\alpha_+(x)<1\quad\mbox{and}\quad\alpha_-:=
\operatorname{\lim\sup}\limits_{|x|\longrightarrow\infty}\alpha_-(x)<1,
\]
and $c_+, c_-\dvtx \R\longrightarrow(0,\infty)$ are such that
\[
\lim_{|x|\longrightarrow\infty}\frac{\alpha_+(x)c_-(x)}{c_+(x)\alpha_-(x)}|x|^{\alpha_+(x)-\alpha_-(x)}=1
\]
and such that condition (\ref{eq14}) is satisfied with the
constant $T(\alpha,\beta)$, where $\alpha:=\alpha_-\vee\alpha_+$ and
$\beta\in(0,1-\alpha)$. In this case, for the test function $V(x)$,
we take $V(x)=1-(1+|x|)^{-\beta}$ again.

In the following proposition, we treat the case when the family of
density functions $\{f_x\dvt x\in\R\}$ is exactly a family of
$S_{\alpha(x)}(\beta(x),\gamma(x),\delta(x))$ densities and we give
sufficient conditions on functions $\alpha(x),\beta(x),\gamma(x)$
and $\delta(x)$ such that the family $\{f_x\dvtx x\in\R\}$ satisfies
conditions (C1$^\prime$)--(C5$^\prime$).
From \cite{taqqu}, Properties 1.2.2, 1.2.3, 1.2.4 and 1.2.15,
\cite{durrett}, Theorem 3.3.5, and (\ref{eq12}) it follows:
%
\begin{proposition}\label{p55} Let $0<\varepsilon<1$, $M>0$ and
$k'\geq
0$ be
arbitrary, and let $F_{\alpha}\subseteq[1,2)$,
$F_\beta\subseteq(-1,1)$ and $F_{\gamma}\subseteq(0,\infty)$ be
arbitrary and finite. Furthermore, let
\begin{enumerate}[(iii)]
\item[(i)] $\bar{\alpha}\dvtx \R\longrightarrow(\varepsilon
,2-\varepsilon)$
and $\tilde{\alpha}\dvtx \R\longrightarrow(0,1)\cup
F_{\alpha}$, such that $\inf_{x\in C}\tilde{\alpha}(x)>0$ for all
compact sets $C\subseteq\R$,
\item[(ii)] $\bar{\beta}\dvtx \R\longrightarrow(-1+\varepsilon
,1-\varepsilon
)$ and $\tilde{\beta}\dvtx \R\longrightarrow
F_{\beta}$,

\item[(iii)] $\bar{\gamma}\dvtx \R\longrightarrow(0,M)$, $\tilde
{\gamma}\dvtx \R
\longrightarrow F_{\gamma}$ and $\hat{\gamma}\dvtx \R\longrightarrow
(\varepsilon,M)$, such
that $\inf_{x\in C}\bar{\gamma}(x)>0$ for all compact sets and
$C\subseteq\R$,
\item[(iv)] $\delta\dvtx \R\longrightarrow(-M,M)$
\end{enumerate}
be arbitrary and Borel measurable. Define
\begin{eqnarray*}
\alpha(x)&:=&\cases{ %
 \bar{\alpha}(x), &\quad $x\in
\bigl[-k',k'\bigr]$,
\cr
\tilde{\alpha}(x), &\quad $x\in \bigl[-k',k'
\bigr]^{c}$,}
\\
\beta(x)&:=&\cases{ %
 \bar{\beta}(x), &\quad $x\in
\bigl[-k',k'\bigr]$,
\cr
\bar{\beta}(x)1_{\{y:\alpha(y)<1\}}(x)+\tilde{\beta}(x)1_{\{
y:\alpha
(y)\geq1\}}(x), &\quad $x\in
\bigl[-k',k'\bigr]^{c}$ }\quad \mbox{and}
\\
\gamma(x)&:=& \cases{ %
 \hat{\gamma}(x), &\quad $x\in
\bigl[-k',k'\bigr]$,
\cr
\bar{\gamma}(x)1_{\{y:\alpha(y)<1\}}(x)+\tilde{\gamma}(x)1_{\{
y:\alpha
(y)\geq1\}}(x), &\quad $x
\in \bigl[-k',k'\bigr]^{c}$. }
\end{eqnarray*}
Then, for any $l'\geq0$, the family of $S_{\alpha(x)}(\beta
(x),\gamma
(x),\delta(x))$, $x\in\R$, densities satisfies
conditions \emph{(C1$^\prime$)--(C5$^\prime$)}.
\end{proposition}
Unfortunately, Proposition \ref{p55} does not cover the case when
the function $\alpha(x)$ takes infinitely many values in the
interval $[1,2)$ since we do not know the series representation
of stable densities for $\alpha\geq1$, as for $\alpha<1$
(see \cite{zolotarev}, Theorems 2.4.2, 2.5.1 and 2.5.4).

At the end, note that all conclusions, methods and proofs
given in this paper can also be carried out in the discrete state
space $\ZZ$. Note that in this case conditions (C1)--(C5) are reduced
just to conditions (C2) and (C3), since compact sets are replaced
by finite sets. Therefore, we deal with a Markov chain
$\{X^{d}_n\}_{n\geq0}$ on $\ZZ$ given by the transition kernel
\[
p_{i,j}=f_i(j-i)
\]
for $i,j\in\ZZ$, where $\{f_i\dvt i\in\ZZ\}$ is a
family of probability functions which satisfies the following
conditions:
\begin{enumerate}[(CD2)]
\item[(CD1)] $f_i(j)\sim c(i)|j|^{-\alpha(i)-1},$
when
$|j|\longrightarrow\infty$, for every $i\in\ZZ$;
\item[(CD2)] there exists $k\in\N$ such that
\[
\lim_{|j|\longrightarrow\infty}\sup_{i\in\{-k,\ldots
,k\}
^{c}}\biggl\llvert f_i(j)
\frac{|j|^{\alpha(i)+1}}{c(i)}-1\biggr\rrvert =0.
\]
\end{enumerate}
Functions $\alpha\dvtx \ZZ\longrightarrow(0,2)$ and
$c\dvtx \ZZ\longrightarrow(0,\infty)$ are arbitrary given functions.
Proofs and assumptions of Theorems \ref{tm13} and \ref{tm14} in
the discrete case remain the same as in the continuous case because
we can switch from sums to integrals due to the tail behavior of
transition jumps.

\section*{Acknowledgements} The author would like to thank Prof. Zoran
Vondra\v{c}ek for many discussions on the topic and for helpful
comments on the presentation of the results. Many thanks to Prof.
Ren\'{e} Schilling and Dr. Bj\"{o}rn B\"{o}ttcher from TU Dresden
for discussions on the topic during their visit to Zagreb and
authors visit to Dresden, funded by DAAD and MZOS of the Republic
of Croatia. The author also thanks the Referees and the Associate
Editor for helpful and constructive criticism of an earlier version
of this paper.



%

\printhistory

\end{document}